\documentclass[12pt]{amsart}
\usepackage[colorlinks,linkcolor=blue,citecolor=blue,urlcolor=blue, pdfcenterwindow, pdfstartview={XYZ null null 1.2}, pdffitwindow, pdfdisplaydoctitle=true]{hyperref}
\usepackage[english]{babel}
\usepackage{array}
\usepackage{threeparttable}
\usepackage{rotating}
\usepackage[section]{placeins}
\usepackage{amssymb} 
\usepackage{mathrsfs}
\usepackage{enumerate}
\usepackage{color}
\usepackage{relsize}
\usepackage{amsrefs}        
\usepackage[all]{xy}

\input xy
\xyoption{all}
\begin{document}
\newtheorem{lemma}{Lemma}[section]
\newtheorem{lemm}[lemma]{Lemma}
\newtheorem{prop}[lemma]{Proposition}
\newtheorem{coro}[lemma]{Corollary}
\newtheorem{theo}[lemma]{Theorem}
\newtheorem{conj}[lemma]{Conjecture}
\newtheorem{prob}{Problem}
\newtheorem{ques}{Question}
\newtheorem{rema}[lemma]{Remark}
\newtheorem{rems}[lemma]{Remarks}
\newtheorem{defi}[lemma]{Definition}
\newtheorem{exam}[lemma]{Example}

\newcommand{\N}{\mathbb N}
\newcommand{\Z}{\mathbb Z}
\newcommand{\R}{\mathbb R}
\newcommand{\Q}{\mathbb Q}
\newcommand{\C}{\mathbb C}
\title{On the bounded cohomology of Lie groups}

\author{I. Chatterji}
\author{G. Mislin}
\author{Ch. Pittet}
\author{L. Saloff-Coste}

\address{MAPMO Universit\'e d'Orl\'eans, Orl\'eans}    
\curraddr{Department of Mathematics, Ohio State University}
\email{Indira.Chatterji@univ-orleans.fr}

\address{Department of Mathematics ETHZ, Z\"urich}    
\curraddr{Department of Mathematics, Ohio State University}
\email{mislin@math.ethz.ch}

\address{CMI Universit\'e d'Aix-Marseille I, Marseille}
\curraddr{Laboratoire Poncelet CNRS, Moscow}
\email{pittet@cmi.univ-mrs.fr}

\address{Department of Mathematics, Cornell University}
\email{lsc@math.cornell.edu}

\keywords{Bounded cohomology, Lie groups, classifying spaces, Borel cocycles, word metrics}

\subjclass[2000]{Primary: 57T10, 55R40; Secondary: 20F65, 53C23}

\thanks{I. Chatterji is partially supported by the NSF grant No. 0644613 and the ANR grant JC-318197 QuantiT, Ch.~Pittet is partially supported by the CNRS. L. Saloff-Coste is partially supported by NSF grant
DMS-0603886.  The authors are grateful to the following institutions were part of this work was achieved: Aix-Marseille I University, the CIRM at Luminy, Cornell University, the Erwin Schr\"odinger Institute, the Institut Henri-Poincar\'e, the Laboratoire Poncelet at the Independent University of Moscow, Ohio State University.}

%
\date{May 9, 2009}

\maketitle
%
\section{Introduction}
\subsection{On  Borel cohomology with $\Z$-coefficients}
Let $G$ be a connected Lie group. 
Recall that the radical $\sqrt G$ of $G$ is its largest connected solvable normal subgroup. Our main result is:

\begin{theo}\label{Main} Let $G$ be a connected Lie group. The following conditions are equivalent.
\begin{enumerate}
	\item\label{RAdLin} The radical $\sqrt G$ of $G$ is linear.
	\item\label{HBb} Each Borel cohomology class of $G$ with   
	            $\Z$-\-coef\-ficients can be represented by a  Borel bounded  
	            cocycle.
	\item\label{H2Bb} Each Borel cohomology class of $G$ of degree two with   
	            $\Z$-\-coef\-ficients can be represented by a  Borel bounded 
	            cocycle.
	\end{enumerate}
\end{theo}

The equivalent conditions of Theorem \ref{Main} admit
geometric (Theorem \ref{theo: main}, Condition \ref{Dis}) and topological
(Theorem \ref{theo: main}, Condition \ref{Com})
counterparts.

Theorem \ref{Main} leads to a generalization (Corollary \ref{coro: generalization of Gromov}) of Gromov's boundedness
theorem \cite{GroVol}*{Section 1.3, p. 23} on characteristic classes of flat bundles.

Before stating our results in more details in Subsection
\ref{subs: main results}, we recall some background
works.

\subsection{Flat bundles and bounded cohomology}
Milnor has shown in \cite{MilOnTheExi} that 
an $SL(2,\mathbb R)$-bundle $P$ over the closed orientable surface $\Sigma_g$ of genus $g>0$
has Euler number satisfying
\begin{equation}\label{Milnor's inequality}
|(e(P),[\Sigma_g])|\leq g-1
\end{equation}
if and only if it admits a flat connection.
Let $G$ be a connected Lie group and let
$c\in H^2(BG,\mathbb Z)$ be a  class of degree $2$
in the cohomology of its classifying space $BG$.
Using the fact that the classifying map of a flat bundle factorizes
through  the classifying space $BG^{\delta}$ of $G$ endowed with the discrete topology, 
Dupont has shown in \cite[Proposition 2.2]{Dup} that for a flat principal $G$-bundle $P$ over $\Sigma_g$ the absolute value of the characteristic number $(c(P),[\Sigma_g])$ is bounded above by the integer part
of $(4g-2)\|c^{\delta}\|_{\infty}$ (where $\|c^{\delta}\|_{\infty}$ is the norm of the image $c^{\delta}\in H^2(BG^{\delta},\mathbb R)$ of $c$):
\begin{equation}\label{Dupont's inequality}
	|(c(P),[\Sigma_g])|\leq \lfloor (4g-2)\|c^{\delta}\|_{\infty}\rfloor.
\end{equation}
In the case $G=SL(2,\mathbb R)$ and $c$ is the Euler class,
$\|c^{\delta}\|_{\infty}=1/4$. Hence Inequality \ref{Dupont's inequality}
and the knowledge of $\|c^{\delta}\|_{\infty}$ imply together Inequality \ref{Milnor's inequality}.  (The upper bound $\|c^{\delta}\|_{\infty}\leq 1/4$ follows from
the bound of the area of hyperbolic triangles, see \cite[Formula 4.3]{Dup}
and \cite[Section 3.1.3]{BucPhd}, or can be directly deduced from
\cite[Lemma 3]{MilOnTheExi}. The lower bound $\|c^{\delta}\|_{\infty}\geq 1/4$, follows from 
Inequality \ref{Dupont's inequality} and the existence, for each $g>1$, of 
a flat $SL(2,\mathbb R)$-bundle over $\Sigma_g$
with Euler number $1-g$, see \cite[Theorem 2]{MilOnTheExi}.) 

The above results raise the following general problem, first studied in degree two by Dupont in \cite{Dup}, then in any degree by Gromov in \cite{GroVol}.
For which Lie groups does the image of the natural map 
$$H^*(BG,\mathbb R)\to H^*(BG^{\delta},\mathbb R),$$
induced by $G^\delta\to G$,
consist of bounded classes? 
For $G$ a connected Lie group, the natural map 
$H^*(BG,\Z)\to H^*(BG^\delta,\Z)$ is always injective, (see Milnor \cite[Cor. 1]{MilOnTheHom}),
contrary to the case of $\R$ coefficients, where it is never
injective, unless $G$ is contractible (it factors through continuous cohomology
$H^*_c(G,\R)$, which is a finite-dimensional $\R$-algebra by van Est's Theorem).
In the case $G_a(\C)$ is the complex Lie group associated to an algebraic group $G_a$ defined over $\mathbb C$, Grothendieck has shown in \cite[Th\'eor\`eme 7.1]{Grot}
that for all fields $k$ of characteristic 
zero, the map 
$$H^d(BG_a(\C),k)\to H^d(BG_a(\C)^{\delta},k),$$
is zero for $d>0$: the
characteristic classes of flat bundles over finite complexes with $G_a(\C)$ as
structure group are torsion classes (see Proposition \ref{prop: flat bundles}). The quotient of the three dimensional complex Heisenberg
group by the subgroup generated by a non-trivial central element shows that
the statement is not true for a general (i.e. not necessary algebraic) connected complex Lie group \cite[Remarques 7.5, b)]{Grot}. In \cite{Gol}, Goldman considers 
the same example ${\bf H}/\Z$ (over $\R$) 
to exhibit non-bounded characteristic classes (in this example, $H^2(B({\bf H}/\Z),\R)\cong\R$  injects into $H^2(B({\bf H}/\Z)^\delta,\R)$). 
Dupont has shown in \cite[Theorem 4.1]{Dup} that 
characteristic classes of degree $2$ of  semi-simple Lie groups with finite center are bounded. He has also given explicit upper bounds in the case of some simple groups (see also \cite{CO} and \cite{COcor}). 
In the case $G$ is the real Lie group associated to a linear algebraic group defined over $\mathbb R$,
Gromov has shown in \cite[Section 1.3, p. 23]{GroVol}, that all characteristic classes are bounded. The result has been improved by Bucher-Karlsson in \cite[Theorem 4]{BucFin}:
each characteristic class of $G$ can be represented by a cocycle whose set of 
values on all singular simplices of $BG^{\delta}$ is finite.

\subsection{Borel cohomology of topological groups}
In \cite[Theorem 10]{Moo}, Moore  has generalized the Eckmann-Eilenberg-Mac Lane theorem to locally compact groups. Namely, he has shown that if $G$ is a locally compact separable group and if $A$ is a polish abelian $G$-module, then the second cohomology group of Borel cochains on $G$ with values in $A$ is isomorphic to the group of  equivalence classes of topological group extensions of $G$ by the $G$-module $A$ (see also Theorem \ref{theo: Moore}): 
\begin{eqnarray*}H^2_B(G,{A})&\cong&{\rm Ext}_{top}(G,{A}).
\end{eqnarray*}
In the case $G$ is a Lie group,  
a result of Wigner (see \cite[Theorem 4]{Wig} and \cite[Section 7]{BSV})
implies that  the cohomology of Borel cocycles on $G$ with values in $\mathbb Z$
(with the trivial $G$-action) is naturally
isomorphic to the singular cohomology of 
the classifying space of $G$ with integer coefficients:
\begin{eqnarray*}H^*_B(G,\mathbb Z)&\cong&H^\ast(BG,\Z).
\end{eqnarray*}

\subsection{Main results}\label{subs: main results}
Building on the works recalled above, we obtain the following theorems.
 
\begin{theo}\label{theo: main} Let $G$ be a connected Lie group. The following conditions are equivalent.
\begin{enumerate}
	\item\label{RadLin} The radical $\sqrt G$ of $G$ is linear.
	\item\label{Com} The closure of the commutator subgroup 
	           of $\sqrt G$ is       
	           simply-con\-nected.
    \item\label{Hb} Each Borel cohomology class of $G$ with $\Z$-\-coef\-ficients 
               can be represented by a  Borel bounded cocycle.
	\item\label{H2b} Each Borel cohomology class of $G$ of degree two with   
	            $\Z$-\-coef\-ficients can be represented by a  Borel bounded cocycle.
	\item\label{Uni} The class in $H^{2}_B(G,\pi_1(G))$ defined by the universal    
	            cover  of $G$ can be represented by a  Borel bounded cocycle.
	\item\label{Dis} The natural inclusion $\pi_1(G)\to\tilde{G}$ of the  
	            fundamental group of $G$ into the universal cover of $G$  
	            is undistorted. 
	
\end{enumerate}
\end{theo}

\begin{rems}
\begin{enumerate}
\item The equivalence between Conditions \ref{RadLin} and \ref{Com} in Theorem \ref{theo: main} is not new: a theorem of Got\^o \cite[Theorem 5]{Got} states that a connected solvable Lie group $S$ is linear if and only if $\pi_1\left(\overline{[S,S]}\right)=0$. 
	\item The proof we shall give, shows that one can relax the boundedness
hypothesis in Conditions \ref{Hb}, \ref{H2b}, and \ref{Uni}, of Theorem \ref{theo: main}: assuming that the representative cocycle has
sub-linear growth, leads in each case to an equivalent condition.
Similarly, assuming that the distortion of $\pi_1(G)\to\tilde{G}$ is sub-linear, 
is equivalent to Condition \ref{Dis} of Theorem \ref{theo: main} (see Definition \ref{defi: distortion}).
	\item Theorem \ref{Main} assumes \emph{integer} coefficients.  Let us emphasis that this is in contrast with the case of \emph{real} coefficients
	by recalling the following well-known points about the forgetful map in Borel bounded cohomology with $\R$-coefficients:
	\begin{enumerate}
	\item $H^*_{Bb}(G,\R)\to H^*_B(G,\R)$ is not onto for $G$ a  connected
	solvable Lie group  for which the right hand side does not vanish for all positive degrees (the left hand side is 0).
	\item $H^*_{Bb}(G,\R)\to H^*_B(G,\R)$
	is not onto in degree 3 for $G$ the
	universal cover of $SL(2,\R)$ (see example \cite[9.3.11 (ii), page 127]{Mon}), thus, in general it is not onto for semi-simple Lie groups either.
	\item $H^*_{Bb}(G,\R)\to H^*_B(G,\R)$
	is conjectured to be onto for semi-simple Lie
	groups with finite center (see \cite[Conjecture 18.1, page 56]{Cha}).
	The conjecture has been proved for Hermitian Lie groups	with finite
	center (see \cite{HartnickOtt}).
	\end{enumerate}
\end{enumerate}
\end{rems}

One ingredient in the proof of Theorem \ref{theo: main} is
the following basic result on the cohomology
algebra $H^*(BG,\R)$ of a connected Lie group.

\begin{theo}\label{theo: higher classes}Let $G$ be a connected Lie group and let $\sqrt G$ be its radical. Then $H^*(BG,\R)$ is generated as an
algebra by $H^2(BG,\R)$ together with the image ${\rm Im}(H^*(B(G/\sqrt G),\R)\to H^*(BG,\R))$. 
\end{theo}

There exist virtually connected (i.e., with finitely many connected components)
Lie groups with non-linear radical and with all classes in $H^2_B(G,\mathbb Z)$ bounded (see Example \ref{counterexample}). Hence Condition \ref{H2b} of Theorem
\ref{theo: main} does
not imply Condition \ref{RadLin} of Theorem \ref{theo: main} for virtually connected Lie groups, but the following holds true.

\begin{theo}\label{theo: vir}
Let $G$ be a virtually connected Lie group. If the radical of $G$ is linear	then each Borel cohomology class of $G$ with $\Z$-coefficients 
can be represented by a Borel bounded cocycle.
\end{theo}

As a corollary,
we obtain the following generalization of Gromov's \cite[Section 1.3, p. 23]{GroVol} and Bucher-Karlsson's \cite[Theorem 4]{BucFin} theorems.

\begin{coro}\label{coro: generalization of Gromov}Let $G$ be a virtually connected Lie group with linear radical. Each class in the image of the natural map $H^*(BG,\R)\to H^*(BG^{\delta},\R)$ 
can be represented by a cocycle whose set of 
values on all singular simplices of $BG^{\delta}$ is finite.
\end{coro}

Let us prove that Corollary \ref{coro: generalization of Gromov} follows from Theorem \ref{theo: vir}.
\begin{proof}
The group $H^*_B(G,\Z)$ is naturally isomorphic to the singular cohomology $H^\ast(BG,\Z)$ of the classifying space $BG$ of $G$ (cf. \cite{Wig} and \cite[Section 7]{BSV}). The following natural diagram commutes,
	$$\xymatrix{H^*_{Bb}(G,\Z)\ar[r]\ar[d]&H^*_B(G,\Z)\cong H^*(BG,\Z)\ar[r]
	   &H^*(BG,\R)\ar[d]\\
	H^*_b(G^\delta,\Z)\ar[r]&H^*(G^\delta,\R)\ar[r]&H^*(BG^\delta,\R),}$$
and the natural isomorphism $H^*(G^\delta,\R)\to H^*(BG^\delta,\R)$ preserves the pseudo-norms (cf. \cite[Section 3.3, p. 49]{GroVol}). Hence, Theorem
\ref{theo: vir} implies that each characteristic class $x^{\delta}$ in the image of the composition
$H^*(BG,\Z)\to H^*(BG,\R)\to H^*(BG^{\delta},\R)$ is bounded. 
Lemma \ref{lemm: bounded implies finite} shows that $x^{\delta}$ can be represented by a cocycle with finite range. 
As $H^*(BG,\R)\cong H^*(BG,\Z)\otimes\R$, we conclude that each characteristic class in the image of $H^*(BG,\R)\to H^*(BG^\delta,\R)$ can be represented by a cocycle with finite range.
\end{proof}

In Corollary \ref{coro: generalization of Gromov}, the hypothesis of linearity on the radical is needed; we will give for every $i>0$ an example of a 
$2i$-dimensional connected Lie group $G_i$ with non-linear radical
such that for all $d$ with $1\le d\le i$,\;
$H^{2d}(BG_i^\delta,\R)$ contains elements in the image
of $H^{2d}(BG_i,\R)\to H^{2d}(BG_i^\delta,\R)$ which cannot be
represented by bounded cocycles (Example \ref{exam: Serre}).

Generalizing to all connected Lie groups a construction of Goldman in \cite{Gol}, we shall prove:

\begin{theo}\label{theo: generalization of Goldman} Let $G$ be a connected Lie group. If the commutator subgroup of the radical of $G$ is not simply-connected relative to its analytic  topology, then there exists a class of degree $2$ in the image of the natural map $H^*(BG,\R)\to H^*(BG^{\delta},\R)$ which can't be represented by a bounded cocycle.  
\end{theo} 

We do not know if the condition 
$$\pi_1\big([\sqrt G,\sqrt G]\big)\neq 0,$$
is the only obstruction to the boundedness of all characteristic classes:
 
\begin{ques}\label{que} Let $G$ be a connected Lie group. Is it true that
	each class in the image of the natural map 
	$H^*(BG,\R)\to H^*(BG^{\delta},\R)$ 
	can be represented by a cocycle whose set of 
	values on all singular simplices of $BG^{\delta}$ is finite, if and only if
	the commutator subgroup of the radical of $G$ is simply-connected relative to  
	its analytic  topology?
\end{ques}

\begin{rema}\label{rema: de Cornulier} \begin{enumerate}
	\item The condition
	$$\pi_1\big([\sqrt G,\sqrt G]\big)\neq 0,$$
	is equivalent to the existence in $\pi_1(G)$ of an  infinite cyclic
	subgroup which is distorted in  $\tilde G$. Starting with a non-trivial element in the fundamental group  of the commutator subgroup of the radical,  the construction in the proof of Proposition \ref{prop: distorted} produces an infinite cyclic distorted subgroup in $\pi_1(G)\subset\tilde G$. The converse implication also holds \cite{DeC}.
	\item 
Let $N$ be a normal analytic subgroup of a connected Lie group $G$.
The inclusion $N\subset\overline{N}$ induces a monomorphism
$$\pi_1\big(N\big)\to\pi_1\left(\overline{N}\right),$$ (see Lemma \ref{closure}). In the case $N=[\sqrt G,\sqrt G]$, if we assume that the monomorphism is an isomorphism, 
Corollary \ref{coro: generalization of Gromov} together with Theorem \ref{theo: generalization of Goldman} bring a positive answer to Question \ref{que}.	See Remark \ref{GuidoEx} for an example of a connected solvable Lie group $G$ with $[G,G]$
simply-connected but $\overline{[G,G]}$ not.
\end{enumerate}
\end{rema}

Question \ref{que} has a positive answer if $G$ is solvable:

\begin{theo}\emph{(Compare with Goldman \cite{Gol}.)}\label{theo: Goldman} Let $G$ be a connected solvable Lie group. The following are equivalent:
\begin{enumerate}
\item\label{commutator} The commutator group $[G,G]$ is simply-connected.
\item\label{virtually trivial} Flat principal $G$-bundles over finite complexes are virtual\-ly trivial.
\item\label{zero} The natural map $H^*(BG,\R)\to H^*(BG^\delta,\R)$ is 0 in all positive degrees.
\item\label{bounded} The image of the natural map $H^*(BG,\R)\to H^*(BG^\delta,\R)$ consist of bounded cohomology classes.
\end{enumerate}
\end{theo}

\subsection{Structure of the paper}
The paper is devoted to proving the theorems stated in Subsection  \ref{subs: main results}. To this end, most of the results recalled and proved in Sections \ref{sect: On the geometry and topology of Lie groups} and \ref{sect: On Borel cohomology}, are needed. Each subsection of these two sections begins with a summary of the main results it contains.
The reader who wishes to understand quickly the main ideas of the paper may
read Section \ref{sect: Proofs of the main results} first, which separately presents the proof of each implication of each theorem of Subsection \ref{subs: main results}, and refer herself/himself to Sections \ref{sect: On the geometry and topology of Lie groups} and \ref{sect: On Borel cohomology}  when needed.

Theorem \ref{theo: main} has roots in \cite[Proposition 5.5 and Lemma 6.3]{CPS}, where some of the authors of the present paper needed to work with Borel cocycles
associated to distorted and undistorted central extensions.
 
\medskip
\noindent
{\it Acknowledgments.} 
We thank Marc Burger, Alessandra Iozzi, Nicolas Monod and Antony Wassermann for useful discussions and comments.
We are grateful to Kenneth Brown for mentioning to us  \cite[Theorem 10]{Moo}.
We are indebted to Yves de Cornulier who explained to
us the geometrical meaning (see Remarks (1) \ref{rema: de Cornulier} and (2) \ref{rema: de Cornulier 2}) of the condition
$$\pi_1\big([\sqrt G,\sqrt G]\big)=0.$$
\section{On the geometry and topology of Lie groups}\label{sect: On the geometry and topology of Lie groups}
\subsection{Preliminaries and facts on Lie groups} In this subsection we discuss general facts about the topology of Lie groups (Lemmata \ref{contractSubgr}, \ref{lemm: pi1} and \ref{lem: universal cover}) and their algebraic structure (Proposition \ref{prop:central}, Lemmata \ref{lem:central}, \ref{lem: one-parameter}, \ref{gensplitting}). We also recall linearity criteria for Lie groups (Theorem \ref{OV} and Proposition \ref{prop: general equivalences}).

To fix notation, we recall some basic facts on Lie groups; a good reference is \cite{OV}. Let G be a connected (real) Lie group. Then G admits a Levi decomposition $G=\sqrt G\cdot L(G)$ with $\sqrt G$ the radical of $G$ (the maximal connected solvable normal subgroup) and $L(G)$ a Levi subgroup (a maximal connected semi-simple subgroup; it is unique up to conjugation). The intersection $\sqrt G\cap L(G)$ is a totally disconnected subgroup of G, which is discrete in $L(G)$ but in general not discrete in $\sqrt G$. In case $G$ is simply-connected, $\sqrt G\cap L(G) = \{e\}$. Hence $G$ is the semi-direct product $G=\sqrt G\rtimes L(G)$. A group $G$ is called \emph{linear}, if it admits a faithful representation $G\to GL(n,\R)$. In case $G$ is a connected  Lie group, there is a closed normal subgroup $\Lambda(G)$, the \emph{linearizer} of $G$, such that $G/\Lambda(G)$ is linear and such that any Lie homomorphism $G\to H$ with $H$ linear factors through $G/\Lambda(G)$ (for the structure
of the linearizer see  \cite{HocTheUni}). If $G$ is connected 
and semi-simple, $\Lambda(G)$ is a central discrete subgroup and any
quotient of $G/\Lambda(G)$ is linear too \cite[Ch. 5, \S 3, Thm. 8]{OV}. In case $G$ is connected and solvable, 
$\Lambda(G)$ is a central torus; a quotient of a linear solvable Lie group need
not be linear. Concerning the linearity of Lie groups, the following theorem is basic.
\begin{theo}\emph{(Malcev \cite{Mal}.)}\label{OV} A connected Lie group $G$ is linear if and only if its radical $\sqrt G$ and Levi subgroup $L(G)$ are.
 \end{theo}
In the sequel we will also deal with not necessarily connected Lie groups $G$; we will write $G^0$ for the connected component of $G$. By the radical $\sqrt G$ of $G$ we still mean a maximal connected normal solvable subgroup of $G$, thus $\sqrt G=\sqrt{G^0}$. A group $G$ is called \emph{virtually connected}, if $G/G^0$ is a finite group. 
In case $G_a$ is a connected linear algebraic group defined over $\R$, its Lie group of $\R$-points $G_a(\R)$ is virtually connected; it might fail to be connected 
as a Lie group. On the other hand every connected linear reductive Lie group $G$ is, as a Lie group, isomorphic to $G_a(\R)^0$ for some connected linear real algebraic group $G_a$ (see \cite{LW}). The following simple observation is used later.

\begin{lemma}\label{contractSubgr} Let $G$ be a contractible Lie group and $H<G$ a connected subgroup. Then $H$ is contractible too.
\end{lemma}

\begin{proof} Recall that any connected Lie group is homotopy
equivalent to a maximal compact subgroup. Let $K$ be a maximal compact subgroup of $H$. Then $K$ is contained in a maximal compact subgroup of $G$, which is $\{e\}$. Thus $K=\{e\}$ and therefore $H$ is contractible.
\end{proof}

Since a maximal compact subgroup of a connected solvable Lie group $S$ is a torus, $S$ is contractible if and only if $S$ is simply-connected. It follows that every connected subgroup of a simply-connected solvable Lie group is contractible.

\begin{lemma}\label{lemm: pi1}Let $N$ be a connected closed normal subgroup of a connected Lie group $G$. Then the inclusion $N\hookrightarrow G$ induces a
short exact sequence $0\rightarrow\pi_1(N)\rightarrow\pi_1(G)\rightarrow
\pi_1(G/N)\rightarrow 0$.\end{lemma}
\begin{proof}The long exact sequence in homotopy reads
$$\cdots\to\pi_2(G/N)\to\pi_1(N)\to\pi_1(G)\to\pi_1(G/N)\to 0\;.$$
Because $G/N$ is a Lie group, $\pi_2(G/N)=0$ (see \cite{Bro}), and
the result follows.\end{proof}
\begin{rema}\label{rema: pi1} Since $\R$ is divisible, ${\rm Hom}(-,\R)$ is an exact contravariant functor on abelian 
groups. Since the fundamental group of a Lie group
is abelian, Lemma \ref{lemm: pi1} implies that the inclusion $N\hookrightarrow G$ induces a surjection ${\rm Hom}(\pi_1(G),\R)\twoheadrightarrow{\rm Hom}(\pi_1(N),\R)$.\end{rema}
\begin{lemma}\label{lem: universal cover}
Let $N$ be a closed connected normal subgroup of a connected Lie group $G$.
There are two exact sequences with commutative squares
$$\xymatrix{
1\ar[r]&\Gamma\ar[r]&\tilde G\ar[r]&G\ar[r]&1\\
1\ar[r]&\Gamma\cap \tilde N\ar[r]\ar[u]&\tilde N\ar[r]\ar[u]&N\ar[r]\ar[u]&1,
}$$
were $\tilde G$ and $\tilde N$ are the respective universal covers of $G$ and $N$, whereas $\Gamma$ and $\Gamma\cap\tilde N$ are central discrete subgroups, isomorphic to the respective fundamental groups of $G$ and $N$, and the vertical arrows are injections of closed normal subgroups.
\end{lemma}	

\begin{proof} Let $p:\tilde G\to G$ be the universal cover of $G$.
The connected component $C$ of the identity in $p^{-1}(N)$ is a closed connected
normal subgroup of $\tilde G$ and the restriction of $p$ to $C$ is a cover
of $N$.	According to Lemma \ref{lemm: pi1} the group $C$ is simply-connected.
Hence $C=\tilde N$, the universal cover of $N$.
\end{proof}

\begin{lemma}\label{closure} Let $p:\tilde G\to G$ be the universal cover of
a connected Lie group $G$. Let $\tilde N$ be an analytic normal subgroup
of $\tilde G$.  Let $N$ denote the projection $p(\tilde N)$ of $\tilde N$ into $G$ with the quotient topology defined by the restriction of $p$ to $\tilde N$. The  inclusion $N\subset \overline N$, induces a monomorphism 
$\pi_1(N)\to \pi_1\left(\overline N\right)$.	
\end{lemma}

\begin{proof} The subgroup $\tilde N$ is closed (see \cite[Chapter 2 (5.14) ]{MM}).
The restriction of $p$ to $\tilde N$ is a cover
of $N$. Hence, according to Lemma \ref{lemm: pi1}, the group $\tilde N$ is the universal cover of $N$. Let $X$ be any path-connected subspace of $G$ containing
$N$ (for example $X=\overline N$). Let us show that the inclusion
$N\subset X$, induces a monomorphism 
$\pi_1(N)\to \pi_1(X)$.
Let $c:S^1\to N$ be a loop such that $c(0)=e$ and assume that
there is a continuous map $h:D^2\to X\subset G$ from the disk to $X$ which extends $c$. Let $\tilde h:D^2\to\tilde G$ be the unique lift of $h$ such that $\tilde h(0)=\tilde e$. The restriction $\tilde c$ of $\tilde h$ to $S^1$ is the unique lift in $\tilde N$ 
of $c$ with $\tilde c(0)=\tilde e$.

\end{proof}

\begin{lemma}\label{lem:central} Let $G$ be a connected nilpotent Lie group.
Let $\phi:\mathbb R\rightarrow G$ be a one-parameter subgroup.
If $\phi(t_0)$ is central for some $t_0\neq 0$, then $\phi(t)$ is central for all $t$.
\end{lemma}

\begin{proof} 
As the universal cover of a Lie group $G$ is an extension of $G$ with discrete central kernel, an easy argument shows that
it is enough to prove the lemma when the nilpotent group $G$ is simply-connected.
In this case, the exponential map
$$\exp:\mathfrak g\rightarrow G$$
is a diffeomorphism. In exponential coordinates the Lie group multiplication is given by the Campbell-Hausdorff
formula	
$$xy=x+y+\frac{1}{2}[x,y]+\frac{1}{12}([x,[x,y]]+[y,[y,x]])+...$$
for all $x,y\in\frak g$. As the Lie algebra $\frak g$ is nilpotent, the above expression has finitely many terms.
Hence if $g\in G$ is given, the equation $g\phi(t)g^{-1}\phi(t)^{-1}=e$ is given by $d=\dim(\frak g)$ polynomial
equations $P_i(t)=0$, $1\leq i\leq d$. By hypothesis, for each $i$, we have  $P_i(nt_0)=0$ for all integer $n$.
Hence the equality holds for all $t\in\mathbb R$.
\end{proof}
\begin{prop}\label{prop:central} A connected nilpotent Lie group $N$ 
has a central torus as its unique maximal compact subgroup.
\end{prop}
\begin{proof} Since a connected compact subgroup of $N$ is a torus, and
since the union of all 1-dimensional tori in a given torus $T$ is
a dense subset of $T$,
it suffices to prove that every 1-dimensional torus $S^1<N$ is central. Choose a surjective homomorphism $\phi: \R\to S^1<N$ and $t_0\in\R\setminus\{0\}$
with $\phi(t_0)=e$. Then the results follows by applying Lemma \ref{lem:central}.
\end{proof}
Recall that the \emph{nilradical} of a Lie group is the largest connected, normal, nilpotent subgroup.
\begin{lemma}\label{lem: one-parameter} Let $K$ be a compact subgroup of a connected solvable Lie group $G$. Let $T_N$ be the maximal compact subgroup of the nilradical $N$ of $G$. If\; $[G,G]\cap K\not=\{e\}$, then there exists a non-trivial one-parameter subgroup in $T_N\cap [G,G]$.\end{lemma}
\begin{proof} Let $k\in [G,G]\cap K$ be a non-trivial element. Let $\frak D$ be the Lie algebra of the (not necessarily closed) analytic subgroup $[G,G]$. Since $G$ is solvable, $[G,G]$ is nilpotent and the exponential map 
$$\exp:{\frak D}\to{[G,G]}$$
is surjective (see \cite{HocTheStr}). Let $X\in\frak D$ such that $\exp(X)=k$ and let $\phi:\R\to{[G,G]}$ be the corresponding one-parameter subgroup $\phi(t)=\exp(tX)$.
Since $[G,G]$ is a normal nilpotent connected subgroup of $G$, it is contained in the nilradical $N$ of $G$. Let $T_N$ be the maximal compact subgroup of $N$. We hence have that ${[G,G]}\cap K\subseteq T_N$. Since $k=\exp(X)\in T_N$, it implies that $\phi(\Z)\subseteq T_N$. The 1-parameter subgroup $\phi(\R)$, being at bounded distance from $\phi(\Z)$ (for any left-invariant Riemannian metric on $G$), is bounded and hence the closure of $\phi(\R)$ in $N$ is a compact subgroup.  We conclude that $\phi(\R)\subseteq T_N\cap [G,G]$.\end{proof}
\begin{prop}\label{prop: general equivalences}Let $G$ be a connected solvable Lie group. The following are equivalent:
\begin{enumerate}\item\label{linc} The group $G$ is linear.
\item\label{commsc} The group $\pi_1\left(\overline{[G,G]}\right)$ is trivial.
\item\label{ccommnc} For every maximal compact subgroup $K<G$, $\overline{[G,G]}\cap K=\{e\}$.
\item\label{commnc} For every maximal compact subgroup $K<G$, $[G,G]\cap K=\{e\}$.
\item\label{Malc} There is a maximal compact subgroup $K<G$ and a closed 1-connected normal subgroup $H<G$ such that $G=HK$. 
\item\label{gMalc} For every maximal compact subgroup $K<G$ there is a 1-connected closed normal subgroup $H<G$ such that $G=HK$ and $H\cap K=\{e\}$.\end{enumerate}\end{prop}
\begin{rema}In  case  the group $G$ is nilpotent, the equivalent conditions of Proposition \ref{prop: general equivalences} are also equivalent to the following one:
{\em{the group is a direct product $G=T\times N$ where $T$ is the unique maximal compact subgroup of $G$ and where $N$ is nilpotent contractible}}.
This follows from Lemma \ref{prop:central}.
\end{rema}

\begin{proof}[Proof of Proposition \ref{prop: general equivalences}] Conditions (\ref{Malc}) and (\ref{gMalc}) are equivalent since any two maximal compact subgroups are conjugate and the equivalence of (\ref{Malc}) and (\ref{linc}) was proved by Malcev in \cite{Mal}. The equivalence of (\ref{linc}) and (\ref{commsc}) was proved by Got\^o, see \cite[Thm.5]{Got}. That (\ref{Malc}) implies (\ref{ccommnc}) and therefore (\ref{commnc}) follows from $\overline{[G,G]}<H$ (to check this inclusion remember that $K$ is abelian), which implies that $\overline{[G,G]}$ is contractible. Since $\overline{[G,G]}$ is homotopy equivalent to its maximal torus, (\ref{ccommnc}) implies (\ref{commsc}) and it remains to show that (\ref{commnc}) implies (\ref{ccommnc}). If (\ref{ccommnc}) does not hold, the maximal compact subgroup $T$ of $\overline{[G,G]}$ is a non-trivial torus with $T\cap [G,G]=\{e\}$ by assumption. Thus there is a continuous embedding $$T\times [G,G]\to\overline{[G,G]},$$ which implies that $\{e\}\times [G,G]$ is dense in 
$T\times [G,G]$, which is a contradiction.
\end{proof}
\begin{rema}\label{GuidoEx} The equivalent conditions of Proposition \ref{prop: general equivalences} are not equivalent to the group $\pi_1\left([G,G]\right)$ being trivial. Indeed, the following is an example of a connected solvable Lie group $G$
with $[G,G]$ simply-connected but $\overline{[G,G]}$ not. Let ${\bf H}$ be the 3-dimensional Heisenberg group and consider ${\bf H}\times S^1$. Its center is $\R\times S^1$. Take the discrete central subgroup $\Z$ generated by $(1,t)$ with 1 generating $\Z$ in $\R$, and $t$ of infinite order in $S^1$; this central subgroup of ${\bf H}\times S^1$ is discrete. Let us define $G:=({\bf H}\times S^1)/\Z$. It is a nilpotent connected Lie group with $[G,G]$ homeomorphic to $\R$, embedded in the maximal 
torus $S^1\times S^1$ of $G$ in a dense way. It follows that $\pi_1([G,G])$ is trivial but $\pi_1\left(\overline{[G,G]}\right)=\Z\times\Z$. Thus $G$ is not linear
(its linearizer is $\overline{[G,G]}$, a central 2-torus). We would like to warn
the reader that in several places in the literature one finds an incorrect
statement saying that a connected solvable Lie group $G$ is linear if and
only if $[G,G]$ is simply-connected (for instance, see \cite[Ch. 2, Thm. 7.1]{OVenc}); the correct statement is
that $G$ is linear if and only if the closure $\overline{[G,G]}<G$ is
simply-connected.
\end{rema}
\begin{lemma}\label{gensplitting} Let $G$ be a connected Lie group and let $\sqrt G$ be its radical. If $Q=G/\sqrt G$ is simply-connected, then the short exact sequence 
$$1\to \sqrt G\to G\to Q\to 1$$
splits, and therefore we have a semi-direct product, $G = \sqrt G\rtimes Q$.\end{lemma}
\begin{proof}This is a classical result in the case $G$ itself is simply-connected (see \cite{HocTheStr}).
The proof in our case is reduced to the simply-connected case by considering the following two exact sequences with commutative squares, 
$$\xymatrix{
1\ar[r]&\sqrt{\tilde{G}}=p^{-1}(\sqrt G)^0\ar[r]\ar[d]&\tilde{G}\ar[r]\ar[d]^p&\tilde{G}/\sqrt{\tilde{G}}\ar[r]
\ar[d]^{\phi}&1\\
1\ar[r]&\sqrt G\ar[r]&G\ar[r]&Q\ar[r]&1,
}$$
where $\tilde{G}$ is the universal cover of $G$. 
The map $\phi$ is an isomorphism because it is a connected covering of a simply-connected space. Thus, if $\sigma$ is a splitting for $\tilde{G}\to
\tilde{G}/\sqrt{\tilde{G}}$, then $p\circ\sigma\circ\phi^{-1}$ is
a splitting for $G\to Q$. \end{proof}

\subsection{Topological $A$-extensions} In this subsection we define {\sl topological $A$-extensions}. Moore's theorem (see Theorem \ref{theo: Moore}) shows that 
they are closely related to $2$-dimensional Borel cohomology.

A central extension of topological groups $0\to A\to E\to G\to 1$ will be referred to as a {\sl topological $A$-extension}. In this setting, we always
assume that the groups are locally compact and second countable and that the monomorphism $A\to E$ has closed image.  
We write ${\rm Ext}_{top}(G,{A})$ for the group of isomorphism classes of topological $A$-extensions
of $G$, where two extensions $0\to{A}\to E_1\to G\to 1$ and $0\to{A}\to E_2\to G\to 1$
are called \emph{isomorphic} if there is an isomorphism of topological groups $\varphi: E_1\to E_2$ making the following diagram commute:
$$\xymatrix{
& &E_1\ar[d]^{\varphi}\ar[dr]& &\\
0\ar[r]&{A}\ar[ur]\ar[r]& E_2\ar[r]&G\ar[r]&1\,.}$$
In case $G$ is a connected Lie group and $A$ is discrete, a topological $A$-extension is just a
covering space $E\to G$ with covering transformation group $A$. Such an $E$ is determined by a
homomorphism $\phi: \pi_1(G)\to A$, with 
$$E=A\times_\phi\tilde{G}=(A\times\tilde{G})/\pi_1(G)$$
with $\tilde{G}$ the universal
cover of $G$, where $\pi_1(G)$ acts on $\tilde{G}$ via deck transformations and on $A$ via $\phi$. It follows that there is a natural isomorphism
$${\rm Ext}_{top}(G,A)\cong {\rm Hom}(\pi_1(G),A)\,.$$

\subsection{Distorted and undistorted central extensions}
In this subsection, we first recall some basic facts about the distortion
of subgroups. The main results are then the following. We prove that in a topological $A$-extension of a Lie group,
the distortion of $A$ is a lower bound for the growth of any Borel cocycle
defining the extension (Proposition \ref{prop: undistorted}).  Then we establish an algebraic criterion to decide when certain
central subgroups of a simply-connected solvable Lie group are distorted (Proposition \ref{prop: Ithaca}). 
Finally, we show that if the radical of a connected Lie group 
$G$ is not linear, then the fundamental group of $G$ is distorted in the universal cover of $G$ (Proposition \ref{prop: distorted}).

\begin{defi}\label{defi: distortion} Let $A$ and $E$ be 
two locally compact, compactly generated groups. We denote by $L_S$, resp. $L_U$, the word length associated to a symmetric relatively compact generating set $S$ of $A$, resp. $U$ of $E$. Assume that $A$ is a subgroup of
$E$. We say that $A$ is \emph{undistorted} in $E$ if the identity map is a quasi-isometry between $(A,L_S)$ and $(A,L_U|_A)$. Otherwise we say that the subgroup $A$ is \emph{distorted} in $E$.\end{defi}
\begin{rems} 
\begin{enumerate}

\item Under the hypothesis of the above definition, Gromov defines in \cite[Chapter 3]{GroAsy}, the distortion function as
$$\hbox{\sc Disto}(r):=\frac{{\rm diam}_A(A\cap B_E(r))}{r}, \forall r>0.$$
One checks that $A$ is undistorted exactly when the function {\sc Disto} is bounded. 
\item We say that the distortion of  $A$ is \emph{at least linear} if there exist $a>0$ and $R>0$ such that for all $r\geq R$, we have $$ar\leq \hbox{\sc Disto}(r).$$
We say that the distortion is \emph{sub-linear} if it is not at least linear.
\item Having a bounded, or unbounded, or sub-linear, etc., distortion function, is a well defined property of the couple
$A<E$, i.e. does not depend on the choice of the relatively compact symmetric generating sets.
\end{enumerate}
\end{rems}
\begin{defi} Let $A$ and $E$ be locally compact, compactly generated groups.
A topological $A$-extension $E$ is called \emph{undistorted}, resp. \emph{distorted}, if $A$ is undistorted, resp. distorted, in $E$. 
\end{defi}

\begin{lemma}\label{disto}Let $p:G\to Q$ be a continuous homomorphism between locally compact, compactly generated groups. Let $H<G$ be a compactly generated subgroup of $G$ and assume that $H\cap\ker(p)=\{e\}$. If $H$ is distorted in $G$, then so is $p(H)$ in $Q$.\end{lemma}
\begin{proof} Let $S$ be a compact symmetric generating set of $H$. Then $p(S)$ is a compact symmetric generating set of $p(H)$. Let $h\in H$. As $p$ is a homomorphism, $L_{p(S)}(p(h))\leq L_S(h)$. As $H\cap\ker(p)=\{e\}$, we also
have  $L_{p(S)}(p(h))\geq L_S(h)$. The proof now follows because
a homomorphism is $C$-Lipschitz with respect to word metrics (for $C\geq 1$ a constant depending on $p$ and on generating sets) and because the projection
of a path between the identity in $G$ and $h\in H$ is still a path between the identity in $Q$ and $p(h)\in p(H)$.
\end{proof}

\begin{lemma}\label{almost locally bounded} Let $0\to A\to E\to G\to 1$, be an exact sequence of second countable locally compact groups, with 
$A$ discrete in $E$. Let $\sigma$ be a Borel section of the projection $p:E\to G$. Let $\mu_E$ be a Haar measure on $E$.
There exists a Borel subset $X$ of $E$ with the following properties:
\begin{enumerate}
	\item $X$ is relatively compact,
	\item $\mu_E(X)>0$,
	\item $X=X^{-1}$,
	\item $\sigma(p(X))$ is relatively compact.
\end{enumerate}
\end{lemma}

\begin{proof} Let $U$ be an open symmetric relatively compact subset of $E$,
small enough such that the covering projection $p$ restricted to $U$, is a homeomorphism onto
its image $V=p(U)$. Let $\mu_E$, resp.   
$\mu_G$, be the Haar measure on $E$, resp. on $G$, such that $\mu_E(U)=1$, resp. such that $\mu_G(V)=1$. The uniqueness of the Haar measure, together with the fact that the homomorphism $p$ is a local homeomorphism, imply that for any Borel subset $A$ of $U$, we have $\mu_E(A)=\mu_G(p(A))$ (see \cite[Int\'egration, chapitre 7, paragraphe 2, num\'ero 7, proposition 10, p. 60]{Bou}). Let $E=\cup_{n\in\N}K_n$, be an exhaustion of $E$ by compact subsets, such that $K_n\subseteq K_{n+1}$.
Let $B_n=\{v\in V: \sigma(v)\in K_n\}=V\cap\sigma^{-1}(K_n)$. It is a Borel
subset of $G$ because $\sigma$ is a Borel map. Notice that
$$\lim_{n\to\infty}\mu_G(V\setminus B_n)=0.$$
We conclude from this, and from the fact that the modular function of $G$ is bounded on $V$, that there exists $N\in\N$, such that $\mu_G(V\setminus B_N)<1/2$, and such that $\mu_G((V\setminus B_N)^{-1})<1/2$.
Let $$Y=V\setminus((V\setminus B_N)\cup(V\setminus B_N)^{-1}).$$
The set $Y$ has the following properties:
\begin{enumerate}
	\item $Y$ is relatively compact,
	\item $\mu_G(Y)>0$,
	\item $Y=Y^{-1}$,
	\item $\sigma(Y)\subseteq K_N$.
\end{enumerate}	
Let $X=p^{-1}(Y)\cap V$. The set $X$ has the required properties.
\end{proof}

\begin{rema} A projection of separable locally compact groups always
admits a locally bounded Borel section \cite[Lemma 2]{Keh}. The point in the above lemma is 
that the given Borel section $\sigma$ is not assumed to be locally bounded. 	
\end{rema}

The following proposition has first been proved by Gersten \cite{Ger}, in
the setting of finitely generated groups.

\begin{prop}\label{prop: undistorted} Let $0\to A\to E\to G\to 1$, be an exact sequence of Lie groups, with $E$ connected,
$A$ finitely generated discrete (hence central) in $E$. Let $\sigma$ be a Borel section of the projection $p:E\to G$, with $\sigma(e)=e$. Let $c:G\times G\to A$, $c(g,g')=\sigma(g)\sigma(g')\sigma(gg')^{-1}$ be the corresponding cocycle. If $c$ is bounded, then $A$ is undistorted in $E$.
\end{prop}

\begin{proof} 
Let $x\in E$. There exist a unique element $a=a(x)\in A$ and a unique element $g=g(x)\in G$ such that $x=a\sigma(g)$.
Let $X\subseteq E$ be as in Lemma \ref{almost locally bounded}.
The set $X(\sigma(p(X)))^{-1}$, is relatively compact in $E$. 
This implies, together with the hypothesis that $A$ is discrete in $E$, 
that there exists a finite symmetric generating set $S$ of $A$
with the property that if $x\in X$, then $a(x)\in S$.

As $\mu_E(X)>0$, the set $XX^{-1}=XX$ contains a neighborhood of the identity
of $E$ (see \cite[Chapter XII, Section 61, Exercise 3, p. 268]{Hal}). As $E$ is connected, the set $X$ generates $E$.

Let $a\in A$ of length relative  to $X$ equal to $n$. Hence, there exist $x_1,\dots, x_n\in X$ such that
\begin{eqnarray*} a&=&x_1\cdots x_n=a_1\sigma(g_1)\cdots a_n\sigma(g_n)\\
&=&a_1\cdots a_n  c(g_1,g_2)\cdots c(g_1\cdots g_{n-1},g_n)\sigma(g_1\cdots g_n).
\end{eqnarray*} 
Applying $p$ to the above equality, we get 
$$e=p(a_1\cdots a_n  c(g_1,g_2)\cdots c(g_1\cdots g_{n-1},g_n)\sigma(g_1\cdots g_n))=g_1\cdots g_n.$$ 
Since $\sigma(e)=e$, we see that the length of $a$ relative to $S$ is bounded
by $n+(n-1)\sup\{L_S(c(g,g'))|L_{p(X)}(g)<n, L_{p(X)}(g')\leq 1\}$.
\end{proof}
\begin{rema}The proof of Proposition \ref{prop: undistorted} shows in fact that
the distortion of $A$ in $E$ is bounded by the growth of $c$.\end{rema}

\begin{prop}\emph{(Gromov \cite[3B2]{GroAsy}.)}\label{prop: linear distortion} Let $N$ be a simply-connected nilpotent Lie group with Lie algebra 
$\frak n$. Let $x\in [\frak n,\frak n]\setminus\{0\}$. Then the one-parameter subgroup $t\mapsto \exp(tx)$, is at least linearly distorted in $N$.
\end{prop}

The statement in \cite[3B2]{GroAsy} which implies Proposition \ref{prop: linear distortion} is not proved.  For a proof of Proposition \ref{prop: linear distortion}, we refer the reader to \cite[Prop. 4.1]{Pit} or \cite{Osi}. A more conceptual proof, in the spirit of \cite{GroAsy}, follows
from the existence of a homothety, relative to a Carnot-Caratheodory metric on $N$, in the case $N$ is graded. The general case is reduced to the graded case
using \cite[Theorem 1.3]{Bre}.   

\begin{prop}\label{prop: Ithaca}Let $0\to\frak n\to\frak g\to{\frak a}\to 0$ be an exact sequence of Lie algebras over $\R$ with $\frak n$ nilpotent and ${\frak a}$ abelian. Let $\frak z\subseteq\frak n$ be a subalgebra which is central in $\frak g$. Let $G$ be the simply-connected Lie group whose Lie algebra is $\frak g$ and let $Z$ be the analytic subgroup of $G$ corresponding to $\frak z$. Then either $\frak z$ is a direct factor in $\frak g$, or there exists a one-parameter subgroup in $Z$ which is distorted in $G$.\end{prop}
\begin{proof}Let $T\in\frak g$ be a regular element, and let $\frak n(T,\frak g)$ be the associated Cartan subalgebra \cite[Chapitre VI, Paragraphe 4, 2, Proposition 9, p. 387]{Che}. The subalgebra $\frak n(T,\frak g)$ is the big kernel of $ad(T)$; that is, there is an integer $N\in\N$, big enough such that $\frak n(T,\frak g)=\ker\left(ad^N(T)\right)=\ker\left(ad^{N+1}(T)\right)$. Let $\frak i(T,\frak g)$ be the small image of $ad(T)$; that is, there is an integer $M\in\N$, big enough such that $\frak i(T,\frak g)=ad^M(T)(\frak g)=ad^{M+1}(T)(\frak g)$. Hence $\frak g$ decomposes as a direct sum of $ad(T)$-invariant subspaces,
$$\frak g=\frak n(T,\frak g)\oplus \frak i(T,\frak g).$$ 
The derivation $ad(T)$ has a semi-simple part $ad_s(T)$, that gives $\frak n(T,\frak g)=\ker\left(ad_s(T)\right)$, and  
$\frak i(T,\frak g)=ad_s(T)(\frak g)$, see \cite[Chapitre VI, Paragraphe 4, 2, Proposition 8, p. 385]{Che}. It follows that
$$[\frak n(T,\frak g),\frak i(T,\frak g)]\subseteq\frak i(T,\frak g),$$
because if $X\in\frak n(T,\frak g)$ and $Y\in\frak i(T,\frak g)$, we can choose $\tilde{Y}\in\frak g$ such that $ad(T)(\tilde Y)=Y$, hence 
$$[X,Y]=[X,ad_s(T)(\tilde{Y})]+[ad_s(T)(X),\tilde{Y}]=ad_s(T)([X,\tilde{Y}]).$$

We note that the Cartan subgroup $C<G$ associated to $\frak n(T,\frak g)$ is,
by Lemma \ref{contractSubgr}, a simply-connected nilpotent group,  because $G$ is a simply-connected solvable Lie group and therefore contractible. If the intersection $\frak z\cap[\frak n(T,\frak g),\frak n(T,\frak g)]$ is non-trivial, then according to Proposition \ref{prop: linear distortion}, there is a one-parameter subgroup in $Z$ which is distorted in the Cartan subgroup $C<G$. This obviously implies
that $Z$ is distorted in $G$ as well, and the proof is finished in this case. 

So we can assume that $\frak z\cap[\frak n(T,\frak g),\frak n(T,\frak g)]=\{0\}$. We notice that as the subalgebra $\frak z$ is central in $\frak g$, it is contained in $\frak n(T,\frak g)$. Hence we can
choose  a complement $V$ for $\frak z\oplus[\frak n(T,\frak g),\frak n(T,\frak g)]$ in $\frak n(T,\frak g)$.
Obviously $\frak m=[\frak n(T,\frak g),\frak n(T,\frak g)]\oplus V$ is a subalgebra of $\frak n(T,\frak g)$ and hence $\frak n(T,\frak g)=\frak z\times\frak m$. Let 
$\frak u=[\frak n,\frak n]\cap\frak n(T,\frak g)$.

We can assume that $\frak z\cap\frak u=\{0\}$. Because otherwise, according to Proposition \ref{prop: linear distortion}, there is a one-parameter subgroup in $Z$ which is distorted in the analytic subgroup of $G$ associated to $\frak n$, hence in $G$ as well. Let $q:\frak z\times\frak m\to\frak m$ be the projection onto the second factor. Let $V'$ be a complement for $\frak u\cap\frak m$ in $\frak u$, that is $(\frak u\cap\frak m)\oplus V'=\frak u$.
If $q(V')\cap([\frak m,\frak m]+\frak u\cap\frak m)\neq\{0\}$, we claim that there is a one-parameter subgroup in $Z$ which is distorted in $G$. To see why, let $v\in V'$ such that $q(v)=x+y\neq 0$, with $x\in [\frak m,\frak m]$ and $y\in\frak u\cap\frak m$. 
The element $z=v-(x+y)$ belongs to $\frak z$. Notice that $z\neq 0$,
because otherwise we would have $v=x+y\in V'\cap(\frak u\cap\frak m)=\{0\}$. Similarly $x\neq 0$, because otherwise $z=v-y\in\frak u\cap\frak z=\{0\}$. 
Since $x\in[\frak m,\frak m]$, the one-parameter subgroup
$t\mapsto\exp(tx)$ is distorted in $G$, according to Proposition \ref{prop: linear distortion}.
Since $v-y\in \frak u\subseteq [\frak n,\frak n]$, the one-parameter subgroup $t\mapsto\exp(t(v-y)x)$ is distorted in $G$, according to Proposition \ref{prop: linear distortion}. The sub-algebra of $\frak g$ spanned by $x$ and $z$ is isomorphic to $\R^2$. Hence, as $z=(-x)+(v-y)$, the one-parameter subgroup $t\mapsto\exp(tz)$ is also distorted in $G$. 
(The geometric picture is the following. To reach the element $\exp(n^2z)$ where
$n$ is large, we start from the identity in $G$ and we first reach with a path  in $\exp([\frak m,\frak m])$, of length linear in $n$, the point $\exp(n^2(-x))$.
Then we follow a path between $\exp(n^2(-x))$ and $\exp(n^2z)$, obtained as the
left-translated in $G$ by $\exp(n^2(-x))$ of a path in $\exp([\frak n,\frak n])$
between the identity and $\exp(n^2(v-y))$ and of length linear in $n$.)
 
Hence we assume that $q(V')\cap([\frak m,\frak m]+\frak u\cap\frak m)=\{0\}$, and we will show that $\frak z$ is a direct factor in $\frak g$. We consider  a complement $W$ for $q(V')\oplus ([\frak m,\frak m]+\frak u\cap\frak m)$ in $\frak m$. Hence
$$\frak m=([\frak m,\frak m]+\frak u\cap\frak m)\oplus q(V')\oplus W.$$
We define the sum of subspaces,
$$\tilde{\frak m}=([\frak m,\frak m]+\frak u\cap\frak m)+V'+W.$$
We claim that $\frak z\times\tilde{\frak m}\cong\frak n(T,\frak g)$.
The inclusion $\frak z+\tilde{\frak m}\subseteq\frak n(T,\frak g)$ is obvious.
The opposite inclusion is true because if  
$x\in\frak n(T,\frak g)$, then $x=z+y$ with $z\in\frak z$ and $y\in\frak m$, and
we can write  $y=y_1+y_2+y_3$ with 
$y_1\in [\frak m,\frak m]+\frak u\cap\frak m$, $y_2\in q(V')$, and $y_3\in W$. 
Let $v\in V'$ such that $q(v)=y_2$. Define $z'=v-y_2\in\frak z$.
Hence $x=z+y=(z-z')+y_1+v+y_3$, with $z-z'\in\frak z$, $y_1\in [\frak m,\frak m]+\frak u\cap\frak m$, $v\in V'$, and $y_3\in W$. The sum is direct because,   
\begin{eqnarray*}\dim(\tilde{\frak m})&\leq &\dim([\frak m,\frak m]+
\frak u\cap\frak m)+\dim(V')+\dim(W)\\
&=&\dim([\frak m,\frak m]+\frak u\cap\frak m)+\dim(q(V'))+\dim(W)=
\dim(\frak m).\end{eqnarray*}
As the subspace $\tilde{\frak m}$ of $\frak n(T,\frak g)$ is a sub-algebra 
(because it contains $[\frak m,\frak m]=[\frak n(T,\frak g),\frak n(T,\frak g)]$),
and as $\frak z$ is central, we obtain a direct product as claimed. What we gained in replacing $\frak m$ with $\tilde{\frak m}$, is that the latest contains $\frak u$ because $\frak u=(\frak u\cap\frak m)\oplus V'\subseteq \tilde{\frak m}$. This will be crucial in finishing the proof.

We have: 
$$\frak g=\frak n(T,\frak g)\oplus \frak i(T,\frak g)=
\frak z\oplus\tilde{\frak m}\oplus \frak i(T,\frak g).$$ 
The proof will be finish if we show that 
$\tilde{\frak m}\oplus \frak i(T,\frak g)$ is a sub-algebra of $\frak g$.
Let $x,x'\in\tilde{\frak m}$ and $y,y'\in\frak i(T,\frak g)$. We have,
$$[x+y,x'+y']=\underbrace{[x,x']}_{\in\tilde{\frak m}}
+\underbrace{[x,y']+[y,x']}_{\in\frak i(T,\frak g)}+[y,y'].$$
Hence we have to show that $[y,y']\in\tilde{\frak m}\oplus \frak i(T,\frak g)$. As $\frak i(T,\frak g)\subseteq\ker(p)=\frak n$, we have $[\frak i(T,\frak g),\frak i(T,\frak g)]\subseteq[\frak n,\frak n]$. As $\frak n$ is an ideal of $\frak g$, it is
preserved by $ad(T)$. Hence, as $ad(T)$ is a derivation it also preserves 
$[\frak n,\frak n]$. Let $\frak n(T,[\frak n,\frak n])$, resp.
$\frak i(T,[\frak n,\frak n])$, be the big kernel, resp. the small image, of the restriction of $ad(T)$ to $[\frak n,\frak n]$. We have
$$[\frak n,\frak n]=\frak n(T,[\frak n,\frak n])\oplus\frak i(T,[\frak n,\frak n]).$$
Now we can conclude because 
$\frak n(T,[\frak n,\frak n])\subseteq\frak n(T,\frak g)\cap[\frak n,\frak n]=\frak u\subseteq\tilde{\frak m}$ and 
$\frak i(T,[\frak n,\frak n])\subseteq\frak i(T,\frak g)$.
\end{proof}

\begin{rema} In fact, the proof of Proposition \ref{prop: Ithaca} shows that if the distortion function of the central one-parameter subgroup under consideration is not bounded, then it grows at least linearly.	
\end{rema}

\begin{lemma}\label{lemm: Euclidean}
Let $V$ be a simply-connected abelian Lie group with a left-invariant Riemannian metric
and let $\Gamma$ be a lattice of $V$. Let $t\mapsto\exp(tX)$ be a one-parameter subgroup of $V$
which projects to a dense subgroup of $V/\Gamma$. Let $\Gamma=\Z\oplus A$ be a splitting of $\Gamma$ and let $z$ be a generator of the $\Z$-factor. Then there are constants $C>1$, $0<\alpha<1<\beta$,
such that for each $n\in\Z$, there exists $t\in\R$, satisfying $\alpha n<|t|<\beta n$ such that in the cylinder $V/A$ equipped with
the Riemannian metric locally isometric to $V$, we have  $d(nz,\exp(tX))\leq C$
(where $nz$, and $\exp(tX)$ are viewed in the cylinder $V/A$).\end{lemma}

\begin{proof} It is enough to prove the lemma in the case $V=\R^d$ with the usual coordinates
and metric, $\Gamma=\Z^d$ generated by the canonical basis vectors, $t\mapsto tX$ with dense image in $\R^d/\Z^d$,
$\Z^d=\Z\oplus\Z^{d-1}$, and $z$ the first vector of the canonical basis.
As $t\mapsto tX$ has dense image in $\R^d/\Z^d$, the vector $X$ is not orthogonal to $z$.  We may assume that 
$X$ and $z$ are in the same connected component of the complement of the hyper-plane
orthogonal to $z$. Let $0\leq\theta<\pi/2$ be the angle between $X$ and $z$. We may assume
$n\in\N$. We choose  $t=z^n/\cos\theta$. In the cylinder $\frac{\R^d}{\{0\}\times\Z^{d-1}}$, we have
$d(nz,tX)\leq \sqrt 2/2$. 
\end{proof}

\begin{prop}\label{prop: distorted} Let $G$ be a connected Lie group and let $K$ be a compact subgroup of the radical $\sqrt G$ of $G$. If $[\sqrt G,\sqrt G]\cap K\not=\{e\}$, then the fundamental group of $G$ is distorted in the universal cover of $G$. Also there exists a distorted  topological ${\Z}$-extension $E$ of $G$.
\end{prop}

\begin{proof} For this proof let us set $R=\sqrt G$. According to Lemma \ref{lem: one-parameter}, there is a non-trivial one-parameter subgroup $\phi:\R\to [R,R]\cap T_N$, where $T_N$ is the maximal compact subgroup of the nilradical $N$ of $G$, which is also the nilradical of $R$.
The closure of $\phi(\R)$ in $T_N$ is a torus $T\subseteq T_N$. Notice that $T_N$ is central
in $G$ because $T_N$ is a normal subgroup of $G$ with discrete automorphism group.
According to Lemma \ref{lem: universal cover}, the sequence of inclusions of closed connected normal subgroups $T\subseteq R\subseteq G$
induces the commutative diagram,
$$\xymatrix{
1\ar[r]&\Gamma\ar[r]&\tilde G\ar[r]&G\ar[r]&1\\
1\ar[r]&\Gamma\cap\tilde R\ar[r]\ar[u]&\tilde R\ar[r]\ar[u]&R\ar[r]\ar[u]&1\\
1\ar[r]&\Gamma\cap \tilde T\ar[r]\ar[u]&\tilde T\ar[r]\ar[u]&T\ar[r]\ar[u]&1,
}$$
where $\tilde G$, $\tilde R$, and $\tilde T$, are resp. the universal covers of $G$, $R$, and $T$, where $\Gamma$, $\Gamma\cap\tilde R$, and $\Gamma\cap\tilde T$, are central discrete subgroups, resp. isomorphic to the fundamental groups of $G$, $R$, and $T$, and where the vertical arrows are injections of closed normal subgroups. Let $\tilde\phi(t)=\exp(tX)$ be the one-parameter subgroup of $\tilde T$ which covers $\phi$. Let $\frak r$ denote the Lie algebra of $\tilde R$. The vector $X\neq 0$ belongs to $[\frak r,\frak r]$ and is central in the Lie algebra of $\tilde G$, hence we can apply Proposition \ref{prop: Ithaca}, with $\frak g=\frak r$, $\frak n=[\frak r,\frak r]$, $\frak a=\frak r/[\frak r,\frak r]$, and $\frak z=\R X$ to deduce that the one-parameter subgroup $\tilde\phi$ is distorted in $\tilde R$.
As the image of $\tilde\phi$ in $\tilde T\cong\R^{dim(T)}$ is a line and as
$\Gamma\cap \tilde T\cong\Z^{dim(T)}$ is cocompact in $\tilde T$, we deduce that $\Gamma\cap \tilde T$ is distorted in $\tilde R$, hence in $\tilde G$ as well. As $\Gamma\cap \tilde T$ is undistorted in $\Gamma$, this shows that $\Gamma$ is distorted in $\tilde G$. 

To prove the existence of a distorted topological $\Z$-extension
of $G$, let us choose a non-zero element $z_0$ in the free abelian group $\Gamma\cap\tilde T$.
It is possible to choose a direct summand $\Z$ in
the finitely generated abelian group $\Gamma$, such that a generator $z$
of the direct summand $\Z$ shares a  non-zero power with $z_0$: there exist $m, n\in\N$, such that $z_0^m=z^n$. This implies that the infinite cyclic groups $(z_0)$ and $(z)$ lie
at bounded distance from each other in $\tilde G$ (with respect to any left-invariant Riemannian metric on $\tilde G$). Let $B$ be a complement of $(z)=\Z$ in $\Gamma$ and let $E=\tilde G/B$. Let $p:\tilde G\to E$ denote the canonical projection.
The connected Lie group $E$ is a topological $\Z$-extension of $G$, with kernel
generated by the image $p(z)\in E$.
The proof will be finished if we prove that the kernel $(p(z))\cong\Z$ is distorted in $E$.
As $\tilde T$ is the universal cover of the closure $T$ of $\phi$,
and as $\tilde \phi$ is  distorted in $\tilde G$, we deduce from Lemma \ref{disto} that the embedding of $p(\tilde \phi)$ in $E$ is  distorted in $E$, and from Lemma \ref{lemm: Euclidean} (applied with $(p(z_0))=\Z$) that $(p(z))=\Z$ itself is  distorted in $E$.
\end{proof}
\begin{rema}\label{rema: de Cornulier 2} 
\begin{enumerate}
	\item 
The above proof shows, under the hypothesis of Proposition \ref{prop: distorted},  that the distortion function of the  fundamental group in the universal cover of $G$ grows at least linearly. Similarly, the kernel of the distorted topological ${\Z}$-extension $E$ of $G$ is at least linearly distorted.
\item In the example of Remark \ref{GuidoEx}, the fundamental group $\pi_1(G)$ is distorted in the universal cover $\tilde G$ of $G$, but each infinite cyclic subgroup of $\pi_1(G)$ is undistorted in $\tilde G$. This last condition implies that the image of the one-parameter subgroup $\phi$ (in the proof of Proposition \ref{prop: distorted}) is not closed.  
\end{enumerate}
\end{rema}

\section{On Borel cohomology}\label{sect: On Borel cohomology}
\subsection{The relationship between the various cohomology groups}\label{subs: relationship}
In this subsection we recall the relationship we need between various cohomology groups. Recall that a map $f:X\to Y$ of topological spaces, is \emph{Borel} if it is measurable with respect to the $\sigma$-algebras generated by the open subsets of $X$ resp. $Y$.
For $G$ any topological group, we write $H^*_B$ resp. $H^*_{Bb}$, $H^*_c$, $H^*_{cb}$, 
$H^*_b$ for its cohomology based on Borel cocycles (resp. cocycles which are Borel bounded, 
continuous, continuous bounded, or just bounded). We refer the reader to
\cite{MooI&II,Moo,MooIV}, \cite{Gui}, and \cite{Mon} for the definition and functorial properties of
these cohomology theories. For the coefficient group $A$ we take a 
metric abelian group; if $A$ is finitely generated abelian, we always assume that the 
metric corresponds to the word metric coming from a finite symmetric generating set. 
All groups will usually be supposed to be separable and locally compact, with topology given by a complete metric 
(occasionally we will also consider non-separable groups like $\R^\delta$).

Our main object of 
study is the forgetful map $H^\ast_{Bb}(G,\Z)\to H^\ast_B(G,\Z)$ for the case of a virtually connected 
Lie group. Note that 
the target group $H^*_B(G,\Z)$ is naturally isomorphic to the singular cohomology $H^\ast(BG,\Z)$ of 
the classifying space $BG$ of $G$ (cf. \cite{Wig} and \cite[Section 7]{BSV}). There are also canonical isomorphisms 
$H^*_c(G,\R)\cong H^*_B(G,\R)$ (see \cite[Theorem 3, p. 91]{Wig} and ), and $H^*_{cb}(G,\R)\cong H^*_{Bb}(G,\R)$ (see \cite[Section 2.3, (2i), p. 15]{BIW} which refers to \cite[Section 4]{Bla}).  We write $G^\delta$ for $G$ considered as
a discrete group. We have a canonical isomorphism $H^*(G^\delta,\R)\cong H^*(BG^\delta,\R)$. The relationship between the various cohomology groups can be expressed by the following natural commutative diagram:
$$\xymatrix{H^*_{Bb}(G,\Z)\ar[r]\ar[dd]&H^*_B(G,\Z)\cong H^*(BG,\Z)\ar[r]&H^*(BG,\R)\ar[d]^{\phi}\\
&H^*_{cb}(G,\R)\cong H^*_{Bb}(G,\R)\ar[r]\ar[d]&H^*_c(G,\R)\cong H^*_B(G,\R)\ar[d]\\
H^*_b(G^\delta,\Z)\ar[r]&H^*_b(G^\delta,\R)\ar[r]&H^*(G^\delta,\R).}$$
The map $\phi$ from the upper right corner is defined as the composition:
\begin{equation}\label{equa: phi}
\phi:H^*(BG,\R)\cong H^*(BG,\Z)\otimes\R\cong H_B^*(G,\Z)\otimes\R\to H_B^*(G,\R).	
\end{equation}

\subsection{On $d$-dimensional Borel cohomology}
In this subsection, we recall and prove several general properties in bounded cohomology
(Lemmata \ref{modn}, \ref{lemm: bounded implies finite}, \ref{lemm: to R}, \ref{hinv}, Corollaries \ref{covers}, \ref{compact}, Proposition \ref{finite index}). In particular, we prove that if $\Gamma$ is a cocompact lattice in a virtually connected Lie group $G$, then a class in the continuous cohomology of
$G$, whose restriction to $\Gamma$ is bounded, admits a continuous bounded representative cocycle (Proposition \ref{prop: continuous representative}).
We then use this fact to deduce from Gromov's theorem \cite[Section 1.3, p. 23]{GroVol}, that
the forgetful map $H^*_{Bb}(G,\Z)\to H^*_B(G,\Z) $, is surjective for $G$ a semi-simple linear Lie group (Proposition \ref{prop: michelle}). We then recall and prove the geometric interpretation of the vanishing of real characteristic classes in positive degrees: if $G$ is a virtually connected Lie group then the natural map $H^*(BG,\mathbb R)\to H^*(BG^{\delta},\mathbb R)$ is zero in positive degree, if and only if all integral primary characteristic classes of flat principal $G$-bundles over finite complexes are torsion classes (Proposition \ref{prop: flat bundles}). Finally we give examples of unbounded characteristic classes
of any (even) degree (Example \ref{exam: Serre}).

\begin{lemma}\label{modn} Let $G$ be a topological group and $x\in H^d_B(G,\Z)$ be such that
for some $n>0$, $nx$ is bounded. Then $x$ is bounded as well.
\end{lemma}
\begin{proof}
The short exact coefficients sequence, $0\to\Z\to\Z\to\Z/n\Z\to 0$, yields long exact cohomology 
sequences, for Borel bounded cohomology as well as for Borel cohomology:
$$\xymatrix{H^{d-1}_{Bb}(G,{\Z/n\Z})\ar[r]\ar[d]& H^d_{Bb}(G,{\Z})\ar[r]^n\ar[d]&H^d_{Bb}(G,{\Z})\ar[r]\ar[d]&
H^d_{Bb}(G,{\Z/n\Z})\ar[d]\\
H^{d-1}_{B}(G,{\Z/n\Z})\ar[r]& H^d_{B}(G,{\Z})\ar[r]^n&H^d_{B}(G,{\Z})
\ar[r]&H^d_{B}(G,{\Z/nZ})
.}$$
Using that the vertical maps are isomorphisms for $\Z/n\Z$-coefficients,
the result follows from a simple diagram chase.
\end{proof}

\begin{lemm}\label{lemm: bounded implies finite}
Let $G$ be a topological group. Let $x^{\delta}$ be a characteristic class of $G$
in the image of the composition
$H^d(BG,\Z)\to H^d(BG,\R)\to H^d(BG^{\delta},\R)$.
If $x^{\delta}$ is bounded then
it admits a representative cocycle whose set of values on
the singular simplices of $BG^{\delta}$ is a finite subset of $\Z$.		
\end{lemm}

\begin{proof} 
The exact sequence of coefficients $0\to\Z\to\R\to S^1\to 0$ admits a 
bounded section. Let us focus on the following two horizontal exact sequences with commutative squares, it induces.
$$\xymatrix{H^d_b(BG^{\delta},{\Z})\ar[r]\ar[d]& H^d_b(BG^{\delta},{\R})\ar[r]\ar[d]&H^d_b(BG^{\delta},S^1)\ar[d]^{=}\\
H^d(BG^{\delta},{\Z})\ar[r]& H^d(BG^{\delta},{\R})\ar[r]&H^d(BG^{\delta},S^1)\,.
}$$
By hypothesis, $x^{\delta}$ has a representative cocycle with integer values, hence
it goes to zero in  $H^d(BG^{\delta},S^1)$. By hypothesis, it is the image
of an element  $w\in H^d_b(BG^{\delta},{\R})$. A simple diagram chase shows that
$w$ is in the image of $H^d_b(BG^{\delta},\Z)\to H^d_b(BG^{\delta},\R)$.
\end{proof}

\begin{lemma}\label{RealCover}
Let $p:H\to G$ be a
covering map of connected Lie groups. Then the induced map $H^*(BG,\R)\to H^*(BH,\R)$
is surjective.
\end{lemma}
\begin{proof} As $p$ is a $\pi_1(G)$-Galois covering and as $\pi_1(G)$ is finitely generated abelian,
we can factor $p$ in a sequence of connected covering spaces $H=X_0\to X_1
\to\cdots X_n=G$ with $X_i\to X_{i+1}$ a cyclic covering. If
$X_i\to X_{i+1}$ is a finite covering, the induced map
$H^*(BX_{i+1},\R)\to H^*(BX_i,\R)$ is an isomorphism, because the
fiber of $BX_i\to BX_{i+1}$ is $\R$-acyclic. In case $X_i\to X_{i+1}$
is an infinite cyclic covering, there is a circle fibration $S^1\to
BX_i\to BX_{i+1}$ with associated Gysin sequence with $\R$-coefficients
$$\xymatrix{H^d(BX_{i+1})\ar[r]^{\cup e}&H^{d+2}(BX_{i+1})\ar[r]&
H^{d+2}(BX_{i})\ar[r]^{\theta(d)}&H^{d+1}(BX_{i+1})}$$
in which $\theta(d)$ is the zero map, because the real cohomology of
a connected Lie group is concentrated in even dimensions. It follows
that $H^*(BX_{i+1},\R)\to H^*(BX_i,\R)$ is surjective and therefore
the composite map $H^*(BG,\R)\to H^*(BH,\R)$ is surjective too.
\end{proof}
This now yields the following useful corollary.
\begin{coro}\label{covers} Let $p:H\to G$ be a
covering map of connected Lie groups. If all Borel cohomology classes in degree $d$ with $\Z$ coefficients for $G$ are
bounded, then they are all bounded for $H$ as well.\end{coro}
\begin{proof} 
Using the natural isomorphism 
$$H^*_B(L,\Z)\cong H^*(BL,\Z)$$
for $L$ a connected Lie group, 
we have a natural commutative diagram
$$\xymatrix{H^d_{Bb}(G,{\Z})\ar[r]^{p^*_{Bb}}\ar[d]& H^d_{Bb}(H,{\Z})\ar[d]\\
H^d_{B}(G,{\Z})=H^d(BG,\Z)\ar[r]^{p^*_{B}}& H^d_{B}(H,{\Z})=H^d(BH,\Z).
}$$
In view of Lemma \ref{RealCover} we know that $p^*_B$ maps
onto a subgroup of finite index. Therefore, if all Borel cohomology 
classes in degree $d$ for $G$ are bounded, then for every
$x\in H^d_B(H,\Z)$ there is an $m>0$ so that $mx$ is bounded. 
But this implies by Lemma \ref{modn}
that already
$x$ is bounded.
  \end{proof}

 \begin{lemma}\label{lemm: to R}
 Let $G$ be an arbitrary topological group. Then $x\in H^*_B(G,\Z)$ is bounded,
 if its image $y\in H^*_B(G,\R)$ in real Borel cohomology
is bounded (i.e., lies in the image of
 $H^*_{Bb}(G,\R)\to H^*_B(G,\R)$). 
 \end{lemma}
 \begin{proof} The short exact sequence $0\to\Z\to \R \to S^1\to 0$ of topological groups admits
 a Borel bounded section $S^1\to \R$ and gives therefore rise to long exact cohomology
 sequences and a commutative diagram
$$\xymatrix{\cdots H^{d-1}_{Bb}(G,S^1)\ar[r]\ar[d]^{=}&H^d_{Bb}(G,{\Z})\ar[r]\ar[d]& H^d_{Bb}(G,{\R})\ar[r]\ar[d]&H^d_{Bb}(G,S^1)\cdots\ar[d]^{=}\\
\cdots H^{d-1}_B(G,S^1)\ar[r]&H^d_{B}(G,{\Z})\ar[r]& H^d_{B}(G,{\R})\ar[r]&H^d_{B}(G,S^1)\cdots\,.
}$$
A simple diagram chase completes the proof.
\end{proof}
\begin{coro}\label{compact} If $G$ is a compact Lie group, then all elements
of $H^*_B(G,\Z)$ are bounded.
\end{coro}
\begin{proof}
Let $x\in H^*_B(G,\Z)$. Then its image $y\in H^*_B(G,\R)$ is obviously bounded,
because by van Est's theorem, the continuous cohomology of a compact Lie
group vanishes in positive dimensions, and $H^*_B(G,\R)=H^*_c(G,\R)$
by Wigner \cite[Theorem 3, p. 91]{Wig}.
\end{proof}

\begin{prop}\label{finite index}
Let $G$ be a virtually connected Lie group and let $G^0$ be the connected component of the identity. If $x\in H^*_B(G,\Z)$
restricts to a bounded class in $H^*_B(G^0,\Z)$ then $x$ itself is bounded. If $y\in H^*(BG^\delta,\R)$ restricts to a class in $H^*(B(G^0)^\delta,\R)$ which has a representing cocycle which takes only finitely many values, then $y$ has such a representative too.
\end{prop}
\begin{proof} Assume $x\in H^*_B(G,\Z)$  restricts to a bounded class 
in $H^*_B(G^0,\Z)$.
It suffices by Lemma \ref{lemm: to R} to show that the image $u=\gamma(x)\in H^*_B(G,\R)$ of $x$ is bounded. We
have a natural commutative diagram with horizontal arrows $\beta$ and $\delta$ induced by
restriction:
$$\xymatrix{H^d_{Bb}(G,\Z)\ar[r]^{\alpha}\ar[d]&H^d_{Bb}(G,{\R})\ar[r]_{\cong\  \ }^{\beta\  \ }\ar[d]^{\phi}& 
H^d_{Bb}(G^0,{\R})^{G/G^0}\ar[d]^{\psi|\rm{inv}}\\
x\in H^d_B(G,\Z)\   \ar[r]^{\  \gamma}&
H^d_{B}(G,{\R})\ar[r]_{\cong\  \ }^{\delta\  \ }& H^d_{B}(G^0,{\R})^{G/G_0}\,.
}$$
That the horizontal arrows $\beta$ and $\delta$ are isomorphisms follows from the Lyndon-Hochschild-Serre
spectral sequences for the short exact sequences $G^0\to G\to G/G^0$, using
that the Borel bounded cohomology and the Borel cohomology with $\R$ coefficients
vanishes for the finite group $G/G^0$ in positive dimensions 
(for the Lyndon-Hochschild Serre spectral sequence
in continuous, resp. continuous bounded, cohomology see \cite{Bla} resp. \cite[Chapter IV, Section 12]{Mon}; 
as we have mentioned at the beginning of Section \ref{sect: On Borel cohomology},  there are natural isomorphisms $H^*_B(G,\R)\cong H^*_c(G,\R)$, resp. $H^*_{Bb}(G,\R)
\cong H^*_{cb}(G,\R)$).
From our assumption it follows that $\delta(u)=v$
is bounded, say $v=\psi(w)$ for some $w\in H^d_{Bb}(G^0,\R)$. By averaging with respect to the $G/G^0$-action we
can form 
$\overline{w}=\frac{1}{|G:G^0|}\sum gw$ where the sum is taken over a set
of coset representatives of $G^0$ in $G$. Then $\overline{w}\in H^d_{Bb}(G^0,\R)^{G/G^0}$
and $\psi(\overline{w})=v$. Thus $\beta^{-1}({\overline{w})}$ is a bounded representative
for $u$ and we are done with the first case.

The case of $H^*(BG^\delta,\R)$ can be dealt with in a similar way, using the fact that the restriction map induces an isomorphism
$H^*(BG^\delta,\R)\to H^*(B(G^0)^\delta,\R)^{G/G^0}\;.$
\end{proof}
\begin{exam}\label{counterexample}
The following is an example of a virtually connected Lie group $G$ with radical $\sqrt G$
non-linear (actually $\sqrt G=G^0$), but all cohomology classes in  $H^2_B(G,\Z)$ bounded and with
an unbounded class in $H^2_B(G^0,\Z)$; thus the converse of Proposition \ref{finite index} does not hold,
and in our Theorem \ref{theo: main} the implication $(2)\Rightarrow (1)$ would not hold
if we assumed $G$ only to be virtually connected. Take ${\bf H}$
to be the three dimensional real Heisenberg group. It admits an involution  $T:{\bf H}\to {\bf H}$ which
induces multiplication by {-1} on the center of ${\bf H}$. In matrix notation,
$$T:\quad \left(\begin{array}{ccc}
1&a&c\\ 0&1&b\\ 0&0&1\end{array}\right)\to
\left(\begin{array}{crr} 1&-a&-c\\ 0&1&b\\ 0&0&1\end{array}\right).$$
Then $T$ preserves
the infinite cyclic central subgroup $\Z$ generated by any chosen non-zero central element. Thus
$T$ induces an involution on ${\bf H}/\Z$, which induces multiplication by ${-1}$ on
$H^2(B({\bf H}/\Z),\Z)\cong \Z$. It follows that the semi-direct product $G:=({\bf H}/\Z)\rtimes_T\Z/2\Z$
has $H^2_B(G,\Z)$ finite, hence bounded according to Lemma \ref{modn} (it is a finite group isomorphic to $H^2(BO(2),\Z)\cong H^2(P^2(\R),\Z)\cong\Z/2\Z$, because $H^2_B(G,\Z)\cong H^2(BG,\Z)$ and $G$ has $O(2)$ as
a maximal compact subgroup). But
it follows from Theorem \ref{theo: main} that
$0\neq x\in H^2_B(G^0,\Z)\cong H^2(BG^0,\Z)\cong \Z$ is
not bounded, because $\sqrt{G^0}=G^0$ and
the closed subgroup $[\sqrt{G^0},\sqrt{G^0}]\cong\R/\Z$ is not simply-connected.
\end{exam}

\begin{lemma}\label{hinv}Let $A$ be a finitely generated abelian group. If $f:G\to Q$ is a homomorphism 
of Lie groups and a homotopy equivalence of the underlying topological spaces, 
then the induced map $f^*:H^*_B(Q,A)\to H^*_B(G,A)$ is an isomorphism.
In particular, if all classes of $H^*_B(Q,A)$ are bounded, then the same is true
for $H^*_B(G,A)$. \end{lemma}
\begin{proof}The map $f$ induces a homotopy equivalence of classifying spaces $Bf:BG\to BQ$.
Indeed, it induces an isomorphism at the level of homotopy groups and the spaces in question
have the homotopy type of $CW$-complexes. On the other hand, $H^*_B(G,A)$ is naturally isomorphic
to $H^*(BG,A)$ (cf. \cite{Wig}). This proves that $f^*$ is an isomorphism.
If $c$ is a bounded representative of $[c]\in H^*_B(Q,A)$, then $f^*c$ is a
bounded representative of $f^*[c]$. This concludes the proof.
\end{proof}

For functorial properties of (bounded) continuous cohomology, we refer the reader to \cite[Chapter III]{Gui} and \cite{Mon}. The following result
also appears in {Bucher-Karlsson \cite[p. 60]{BucPhd}.}
\begin{prop}\label{prop: continuous representative}
Let $G$ be a locally compact group with  a cocompact lattice
$\Gamma$. Let $res:H^*_c(G,\R)\to H^*(\Gamma,\R)$ be the restriction map
and let  $x\in H^*_c(G,\R)$. If $res(x)$ is bounded then $x$ admits a 
continuous, bounded cocycle representative.
\end{prop}

\begin{proof} 
There is a commutative diagram:
$$\xymatrix{
H^*_{cb}(G,\R)\ar[r]^F&H^*_c(G,\R)
\\
H^*_b(\Gamma,\R)\ar[u]^{tr}\ar[r]^F&H^*(\Gamma,\R)\ar[u]^{tr},}$$
where $F$ denotes forgetful maps from bounded cohomologies and $tr$ denotes the transfer maps (cf. \cite[p. 107]{Mon}). The equality
$$tr\circ res=id,$$ holds on $H^*_c(G,\R)$ (as well as on $H^*_{cb}(G,\R)$).
By hypothesis, there exists $y\in H^*_b(\Gamma,\R)$ such that $F(y)=res(x)$. A continuous, bounded cocycle representative of the transfer of $y$ exists by definition and
$$F(tr(y))=tr(F(y))=tr(res(x))=x.$$
\end{proof}

Recall the map $\phi:H^*(BG,\R)\to H^*_c(G,\R)$ from Subsection \ref{subs: relationship}, Equation~\ref{equa: phi}. 

\begin{prop}\emph{(Compare with Bucher-Karlsson \cite[p. 60]{BucPhd}.)}\label{prop: michelle} Let $G$ be a Lie group of the form $G_a(\R)^0$
for some semi-simple linear algebraic $\R$-group $G_a$.
\begin{enumerate}
	\item Every element
in the image of $\phi:H^*(BG,\R)\to H^*_c(G,\R)$ admits
a continuous, bounded cocycle representative.
    \item The forgetful map $H^*_{Bb}(G,\Z)\to H^*_B(G,\Z)$ is surjective.
\end{enumerate}
\end{prop}
\begin{proof} Let $x\in H^*(BG,\R)$. Its image $\phi(x)$  maps to the element  $x^{\delta}\in H^*(BG^\delta,\R)\cong H^*(G^{\delta},\R)$. The proof of Gromov's result \cite[Section 1.3, p. 23]{GroVol} by Bucher-Karlsson \cite{BucFin} works in a semi-algebraic setting, hence applies
to $G=G_a(\R)^0$.  Therefore, the class 
$x^{\delta}$ is bounded. According to 
\cite{BorCom}, there exists a
cocompact lattice $\Gamma$ in $G$. As the arrows of the following commutative diagram,
$$\xymatrix{
H^*_{c}(G,\R)\ar[d]\ar[dr]^{res}\\
H^*_b(G^{\delta},\R)\ar[r]^{res}&H^*(\Gamma,\R),}$$
preserve boundedness, the class $res(\phi(x))\in H^*(\Gamma,\R)$ is also bounded. Proposition \ref{prop: continuous representative} shows that the class $\phi(x)$ 
has a continuous, bounded representative. This proves the first statement
of the proposition.

To prove the second statement, let $x\in H^*_B(G,\Z)$. According to Lemma
\ref{lemm: to R}, it is enough to show that the image $y\in H^*_B(G,\R)$ of
$x$ is bounded. In view of the commutative diagram  of 
Subsection \ref{subs: relationship}, and the first statement
of the proposition, the class $y$ admits a continuous, bounded
cocycle representative.  
\end{proof}

\begin{lemm}\label{inductive}
Let $W$ be a $CW$-complex. If a class $x\in H^d(W,\R)$ in the real singular
cohomology of $W$ is non-zero, then there exists a finite subcomplex $F\subset W$
such that the restriction of $x$ to $F$ is a non-zero element of $H^d(F,\R)$.	
\end{lemm}

\begin{rema} The corresponding statement with $\R$ replaced by $\Z$ is   not true \cite[p. 1212]{McG}.
\end{rema}

\begin{proof} We prove the lemma. We have:
	\begin{eqnarray*} H^d(W,\R)&\cong&Hom_{\R}(H_d(W,\R),\R)
	\cong Hom_{\R}(\lim_{\rightarrow}H_d(F_{\alpha},\R),\R)\\
	&\cong&\lim_{\leftarrow}Hom_{\R}(H_d(F_{\alpha},\R),\R)
	 \cong\lim_{\leftarrow}H^d(F_{\alpha},\R),
	\end{eqnarray*}
where the inductive, resp. projective, limit is taken over the inductive system of finite subcomplexes $F_{\alpha}$ of $W$, resp. over the projective system $Hom_{\R}(H_d(F_{\alpha},\R))$.

\end{proof}

\begin{prop}\label{prop: flat bundles} Let $G$ be a virtually connected Lie group.
	\begin{enumerate}
		\item If the natural map $H^*(BG,\R)\to H^*(BG^{\delta},\R)$
		is zero in positive degree, then all integral characteristic
		classes of  flat principal $G$-bundles over  finite complexes are torsion classes.
		\item\label{restriction is torsion} If $c(EG)\in H^*(BG,\Z)$ is a characteristic class with the property  that  for any flat principal $G$-bundle $P$ over any finite complex $X$, the class $c(P)\in H^*(X,\Z)$ is a torsion class, then $c(EG)$ is in the kernel of the composition
		$$H^*(BG,\Z)\to H^*(BG,\R)\to H^*(BG^{\delta},\R).$$
	\end{enumerate}
\end{prop} 
\begin{proof} Let $P$ be a flat principal $G$-bundle over a finite complex $X$.
We have a commutative diagram: 
$$\xymatrix{
P\ar[d]\ar[r]&EG^{\delta}\ar[d]\ar[r]&EG\ar[d]\\
X\ar[r]&BG^{\delta}\ar[r]&BG,}$$
where vertical arrows are projections of principal $G$-bundles.
It induces the commutative diagram:
$$\xymatrix{
H^*(BG,\Z)\ar[d]\ar[r]&H^*(BG^{\delta},\Z)\ar[d]
\ar[r]&H^*(X,\Z)\ar[d]\\
H^*(BG,\R)\ar[r]&H^*(BG^{\delta},\R)
\ar[r]&H^*(X,\R),}$$
in which characteristic classes are preserved. As $X$ is a finite complex,
the kernel of the vertical arrow on the right is the torsion
subgroup. 

To prove the first claim of the proposition,
let $c(P)\in H^*(X,\Z)$ be a characteristic class. By definition it is the image of $c(EG)\in H^*(BG,\Z)$.
The commutativity of the diagram and the hypothesis that the map $H^*(BG,\R)\to H^*(BG^{\delta},\R)$ is zero imply together that $c(P)$ is a torsion class.
This proves the first part of the proposition.

To prove the second claim of the proposition, we consider $c(EG)^{\delta}\in H^*(BG^{\delta},\Z)$ and 
assume that its image in $H^*(BG^{\delta},\R)$ is not zero. According to Lemma
\ref{inductive}, there exists a finite sub-complex $X\subset BG^{\delta}$, such that our characteristic class $c(P)\in H^*(X,\Z)$ defined by the flat principal $G$-bundle induced by the composition $X\subset BG^{\delta}\to BG$ is not a torsion class. This concludes the proof of the proposition.  
\end{proof}

We recall the following fact about transfer in cohomology between covering spaces.
\begin{lemma}\label{lemm: transfer for covers}
Let $Y\to X$ be a $d$-fold covering space of $CW$-complexes. There is a transfer map:
\[tr: H^*(Y ,\Z)\to H^*(X,\Z),\]
such that the composition
\[H^*(X,\Z)\to H^*(Y,\Z)\to H^*(X,\Z),\] 
is multiplication by d.
Thus, if $x\in H^*(X,\Z)$ maps to $0$ in $H^*(Y,\Z)$, then $d\cdot x=0$.
\end{lemma}
\begin{proof}
Choose a $CW$-structure on $X$ and lift it to $Y$ so that each $n$-cell in $X$ has precisely $d$ counter images in Y. Define $tr$ on the cellular cochain level by mapping a cellular $n$-cochain $f$ on $Y$ to the cellular $n$-cochain $g$ on $X$ defined by 
\[g(\sigma)=\sum f(\sigma_i),\] where $\sigma _i$ runs over the $d$ counter images of $\sigma$. The rest is obvious.		
\end{proof}

The following family of examples of Lie groups with non-linear
radical shows that one can have unbounded characteristic
classes in any positive (even) degree. (The basic idea in the construction
appears in \cite{Gol} and also in \cite[Remarques 7.5, b)]{Grot}.)
 
\begin{exam}\label{exam: Serre} Let $G_i$ be the $i$-fold cartesian product of Heisenberg quotients $\bf{H}/\Z$, where $\bf{H}$ is the 3-dimensional Heisenberg group and $\Z<\bf{H}$
a central subgroup. Note that the maximal compact subgroup of $\bf{H}/\Z$ is $S^1$ so that $H^*(BG_i,\R)$ is a polynomial algebra generated by $2$-dimensional classes
$x_j$, $1\le j\le i$, with the property that the map induced by the $j$-th injection 
${\bf H}/\Z\to G_i$ maps
$x_j$ to a non-zero element in $H^2(B(\bf{H}/\Z),\R)$. According to Goldman 
\cite{Gol} there is a flat ${\bf{H}}/\Z$-bundle over the $2$-torus 
$T^2$ given by a homomorphism
$\phi:\Z\times\Z\to {\bf{H}}/\Z$ such that $\phi^*: 
H^2(B({\bf{H}}/\Z),\R)\to H^2(T^2,\R)$ is
non-zero. Taking $i$-fold products, we get a flat $G_i$-bundle of the $2i$-torus
$T^{2i}$, with associated map $H^*(BG_i,\R)\to H^*(T^{2i},\R)$ mapping the
product $y_d:=x_1\cdots x_d$ non-trivially into $H^{2d}(T^{2i},\R)$,
where $1\le d\le i$. It follows
that the images $y^{\delta}_d\in H^{2d}(BG_i^\delta,\R)$ of $y_d$\;
are all non-trivial. But they cannot have bounded cocycle
representatives, because the bounded cohomology of the nilpotent group
$G_i^\delta$ vanishes in positive dimensions.
\end{exam}

\subsection{On $2$-dimensional Borel cohomology}\label{subs: bounded2cocycles}
In this subsection we first recall Moore's theorem (Theorem \ref{theo: Moore}) which generalizes the Eck\-mann-\-Eilenberg-\-Mac Lane theorem to locally compact groups. In Lemmata \ref{lemm: cocycle defined by a section}, \ref{lemm: pullback}, \ref{lemm: composition}, 
we establish some properties of $2$-cocycles associated to central extensions.
In Proposition \ref{prop: universal cover} we show that the $2$-dimensional Borel cohomology
with $\Z$-coefficients of a Lie group $G$ is bounded if and only if the class defined by  the universal cover of $G$ is bounded. Lemma \ref{lemm: Kunneth}
is a $2$-dimensional K\"unneth formula in Borel bounded cohomology with $\Z$-coefficients. Proposition \ref{prop: onto} shows that for a connected normal closed subgroup $N$ of a connected Lie group $G$, each real class of degree $2$
on the classifying space $BN$ is the restriction of a class on $BG$. The subsection ends with a proof
based on Moore's theorem, of the surjectivity of the forgetful map
$H^2_{Bb}(G,\Z)\to H^2_B(G,\Z)$, for $G$ a real linear algebraic group.

A topological $A$-extension $0\to{A}\to E\to G\to 1$ always admits a Borel  section $\sigma$ (cf. \cite{Keh}).  We denote $c_\sigma: G\times G\to A$ the Borel 2-cocycle 
$$c_{\sigma}(g,g')=\sigma(g)\sigma(g')\sigma(gg')^{-1}.$$
The following statement is a special case of \cite[Theorem 10]{Moo}:
\begin{theo}\emph{(Moore \cite[Theorem 10]{Moo}.)}\label{theo: Moore} Let $G$ and $A$ be locally compact separable groups with $A$ abelian.
The map
\begin{eqnarray*}{\rm Ext}_{top}(G,{A})&\to&H^2_{B}(G,{A})\\
\{0\to{A}\to E\to G\to 1\}&\mapsto&[c_{\sigma}],\end{eqnarray*}
is an isomorphism.\end{theo}

\begin{lemma}\label{lemm: cocycle defined by a section}
Let $G$ be a topological group and let $A$ be a topological abelian group. Let $c:G\times G\to A$ be an inhomogeneous Borel $2$-cocycle  (that is a Borel map such that,
for all $x, y, z\in G$,
$$dc(x,y,z)=c(y,z)-c(xy,z)+c(x,yz)-c(x,y)=0).$$
Assume that $c$ is cohomologous to $c_{\tau}$, where $\tau$ is a Borel
section of a topological $A$-extension $E$ of $G$ (if $G$ and $A$ are locally compact and separable, such an extension with such a section always exist according to Theorem \ref{theo: Moore}). Then there is a Borel section $\sigma$ of $E$, such that $c=c_{\sigma}$.
\end{lemma}

\begin{proof}
By hypothesis, there exists
a Borel map $b:G\to A$, such that $c=c_{\tau}+db$.
Recall that the coboundary operator on a (inhomogeneous) 1-cochain $b$ 
(and with $A$ a the trivial $G$-module) is
$db(g,g')=b(g')-b(gg')+b(g)$. 
The Borel cocycle $c_{\sigma}$, associated to the Borel section $\sigma:G\to E$, defined as $\sigma(g)=\tau(g)b(g)$, is equal to $c$. 
\end{proof}
\begin{rema} Topological $\R$-extensions are topologically
split and can therefore be described by a continuous cocycle $G\times G\to \R$.
\end{rema}

\begin{lemma}\label{lemm: pullback} Let $p:G\to Q$ and $\pi:\tilde{Q}\to Q$ be surjective continuous homomorphisms between Hausdorff topological groups. 
(We view $\pi:\tilde{Q}\to Q$ as a principal $\ker \pi$-bundle.)
\begin{enumerate}
\item\label{pullback exact sequence}	 The total space of the pull-back $p^*(\tilde{Q})\rightarrow G$ of $\tilde{Q}\to Q$ is a closed subgroup
of $G\times\tilde{Q}$ and there is a short exact sequence of topological groups
$$1\to\ker(p)\to p^*(\tilde{Q})\to\tilde{Q}\to 1.$$
\item\label{pullback section} Assume that both $\ker \pi$ in $\tilde{Q}$ and $ker(p)$ in $G$
are central. If the extension
$$0\to\ker(\pi)\to \tilde{Q}\to Q\to 1$$
can be defined by a Borel bounded cocycle, then the same is true for
the extension
$$0\to\ker(\pi)\to p^*(\tilde{Q})\to G\to 1.$$
\end{enumerate}
\end{lemma}

\begin{proof} To prove (1), recall that by definition, $$p^*(\tilde{Q})=\{(g,\tilde{q})\in G\times\tilde{Q}\hbox{ such that }p(g)=\pi(\tilde{q})\}.$$ We embed $ker(p)$ in $p^*(\tilde{Q})$ by composing
	the inclusions $ker(p)\subset G\subset G\times\{e\}$. The map $\tilde{p}:p^*(\tilde{Q})\to\tilde{Q}$, defined as  $(g,\tilde{q})\mapsto \tilde{q}$, gives the wanted exact sequence. 
To prove (2), we apply Lemma \ref{lemm: cocycle defined by a section} in order
to obtain a Borel section $\sigma:Q\to\tilde{Q}$ of $\pi$
such that the associated cocycle $c_{\sigma}$ is bounded. Then $\sigma^*:G\to p^*(\tilde{Q})$ defined as
$g\mapsto(g,\sigma p(g))$, is a Borel section of the projection $p^*(\tilde{Q})\rightarrow G$ and its associated cocycle satisfies
$$c_{\sigma^*}(g,g')=(e,c_{\sigma}(p(g),p(g'))).$$
\end{proof}

\begin{lemma}\label{lemm: composition} Let $p:X\to Y$ and $q:Y\to Z$ be two surjective
homomorphisms of  groups with  central kernels. Let $\sigma$ be a
section of $p$, and
$\tau$ one for $q$. Assume that the associated cocycles $c_\sigma$ and $c_\tau$
have finite range. If the image under $\sigma$ of the kernel of $q$ is central in $X$, then the cocycle
$c_{\sigma\tau}$ associated to the section $\sigma\tau$ of $qp$  has finite range.\end{lemma}
\begin{proof} Since $\sigma\left(\ker(q)\right)$ is central, the kernel of $qp$ is central so that $qp: X\to Z$
is a central $\ker(qp)$ extension. One 
verifies that
$$c_{\sigma\tau}(x,y)=c_\sigma(\tau(x),\tau(y))\cdot \sigma(c_\tau(x,y))\cdot [c_\sigma(c_\tau(x,y),\tau(xy))]^{-1}.$$
Thus $c_{\sigma\tau}$ is a product of three functions each of which takes only finitely many values.
It follows that $c_{\sigma\tau}$ has finite range.
\end{proof}

\begin{rema}\label{facile} 
Consider the following commutative diagram of groups 
$$\xymatrix{
G\ar[dr]^q\ar[r]^p&E\ar[d]^r&\\
 &Q,}$$
with $q$ and $r$ surjective. If $\sigma$ is a section of $q$ then $p\sigma$ is a section of $r$ and 
\begin{eqnarray*}c_{p\sigma}(x,y)&=&p\sigma(x)p\sigma(y)\left(p\sigma(xy)\right)^{-1}\\
&=&p\left(\sigma(x)\sigma(y)\left(\sigma(xy)\right)^{-1}\right)=(p\circ c_{\sigma})(x,y).\end{eqnarray*}
\end{rema}

\begin{prop}\label{prop: universal cover} Let  $G$ be connected Lie group. The following conditions are equivalent.
\begin{enumerate}
\item The class in $H^2_B(G,\pi_1(G))$, corresponding (in the sense of Theorem    
      \ref{theo: Moore}) to the  universal cover of $G$, is bounded. 
\item All classes in $H^2_B(G,\Z)$ are bounded.
\item For every finitely generated abelian group $A$, all classes in \\ $H^2_B(G,A)$  are bounded.
\end{enumerate}
\end{prop}

\begin{proof} We first prove that (1) implies (3).
Let $x\in H^2_B(G,A)$ and let $E$ be a central extension of $G$ with central discrete kernel $A$ corresponding to $x$. Let $E^0$ be the connected component of $E$. The restriction of the projection $E\to G$ to $E^0$ is a connected cover of $G$. Hence, we obtain a commutative diagram of covering groups,
$$\xymatrix{
\tilde{G}\ar[r]\ar[dr]&E\ar[d] \\
&G,}$$
where $\tilde{G}$ is the universal cover of $G$.
Applying the hypothesis together with Remark \ref{facile}, we deduce the existence
of a Borel section $\tau:G\to E^0\subseteq E$ of the projection $E\to G$ such that the corresponding cocycle $c_{\tau}:G\times G\to A$ is bounded. As $[c_{\tau}]=x$,
we conclude that $x$ is bounded.

The fact that (3) implies (2) is trivial and obviously (3) implies (1), since $\pi_1(G)$
is a finitely generated abelian group. It remains to see that (2) implies (3).  Let
$A=\Z^n\times F$ be a finitely generated abelian group,  $F$ a finite group, and let $c:G\times G\to A$
 be a 2-cocycle; it has components $c_i: G\times G\to \Z$, for $1\le i\le n$ and a component $c_F$, taking values in $F$.  The latter is obviously bounded. By assumption, the cocycles $c_i$ are cohomologous to bounded cocycles $d_i$. Thus, the cocycle $G\times G\to A$ with components $d_i$ and $c_F$ is bounded and cohomologous to $c$, finishing the proof.
 \end{proof}

\begin{lemma}\label{lemm: Kunneth}If $G$ and $H$ are two virtually connected Lie groups such that 
the maps $H^2_{Bb}(G,\Z)\to H^2_{B}(G,\Z)$ and $H^2_{Bb}(H,\Z)\to H^2_{B}(H,\Z)$ are onto, then 
so is the map $H^2_{Bb}(G\times H,\Z)\to H^2_{B}(G\times H,\Z)$.\end{lemma}
\begin{proof} Because $\pi_1(BG)=\pi_0(G)$ is finite, $H^1(BG,\Z)=0$, and similarly $H^1(BH,\Z)=0$.
By the K\"unneth formula in singular cohomology, we conclude that the natural map 
$BG\sqcup BH\to BG\times BH= B(G\times H)$ induces an isomorphism $\psi$ at the level of degree 2
cohomology with $\Z$ coefficients. On the other hand, the inclusion $i_G: G\to G\times H$
induces a surjection $i_G^*$ in Borel bounded cohomology (because there is a retraction
$G\times H\to G$), and similarly for the inclusion $i_H$. The following diagram commutes because
for each component of the vertical maps we have a natural commutative diagram:
$$\xymatrix{
H^2_{Bb}(G\times H,\Z)\ar[r]\ar[d]^{(i_G^*,i_H^*)}_{\rm epi}& H^2_B(G\times H,\Z)
=H^2(B(G\times H),\Z)\ar[d]_{\rm iso}^{\psi}\\
H^2_{Bb}(G,\Z)\oplus H^2_{Bb}(H,\Z)\ar[r]^{\rm epi}& H^2(BG,\Z)\oplus H^2(BH,\Z).}$$
This shows that the map $H^2_{Bb}(G\times H,\Z)\to H^2_B(G\times H,\Z)$ is onto as well.
\end{proof}

\begin{prop}\label{prop: onto}
Let $G$ be a connected Lie group and let $N$ be a closed connected normal subgroup of $G$. The natural map
$$H^2(BG,\R)\to H^2(BN,\R),$$
induced by the inclusion $N\subset G$, is onto.	
\end{prop}
\begin{proof}
There is a commutative diagram:
$$\xymatrix{
H^2(BG,\Z)\ar[r]\ar[d]_{\cong}& H^2(BN,\Z)\ar[d]_{\cong}\\
Hom(\pi_1(G),\Z)\ar[r]& Hom(\pi_1(N),\Z)
.}$$
Tensoring with $\R$, we see that the proof will be finished if we show that
$$Hom(\pi_1(G),\R)\to Hom(\pi_1(N),\R)$$
is onto. But this is the case according to Remark \ref{rema: pi1}.
\end{proof}

\begin{theo}\label{theo: bounded extensions}Let ${G_a}$ be a linear algebraic group defined
over ${\mathbb R}$. Let ${K_a}$ and ${S_a}$ be algebraic subgroups of
${G_a}$ defined over ${\mathbb R}$ with the following properties.
\begin{enumerate}
\item The Lie group $K_a(\R)$ is compact.
\item The Lie group $S_a(\R)$ is simply-connected.
\item $G_a(\R)=S_a(\R)K_a(\R)$.
\item $S_a(\R)\cap K_a(\R)=\{e\}$.
\end{enumerate} 
Let  $A$ be a finitely generated abelian group.
Let $G$ be a virtually connected Lie group.
Assume that the connected Lie groups  $G^0$ and $G_a(\R)^0$ are  isomorphic. Then the natural map,
$$H_{Bb}^2(G,A)\to H_{B}^2(G,A),$$
is surjective.
\end{theo}
\begin{proof} Thanks to Lemma \ref{finite index}, we see that we can assume that
$G$ is connected. In this case $G=G^0=G_a(\R)^0$. According to Proposition \ref{prop: universal cover}, it is enough to prove that the universal cover
$p:\tilde G\to G$ of $G$ admits a Borel section whose associated cocycle is bounded. Let $K=K_a(\R)^0$ and $S=S_a(\R)=S_a(\R)^0$. Let
$\tilde K=p^{-1}(K)$. Since $A:=\pi_1(G)\subseteq \tilde K$,
the regular $A$-cover
$$0\to A\to{\tilde G}\to G\to 1,$$
restricts to
$$0\to A\to{\tilde K}\to K\to 1.$$
Let $\tilde S$ be the connected component of the identity
of $p^{-1}(S)$. The subgroup $\tilde S\subset\tilde G$ is closed.
Hence it is a Lie subgroup. Since $S$ is simply-connected,
the restriction of the covering map $p$ to $\tilde S$ is
an isomorphism between $\tilde S$ and $S$. In what follows, if 
$s\in S$, we will denote
by $\tilde s$ the unique element of $\tilde S$ such that $p(\tilde s)=s$.
From this isomorphism and from the hypothesis
$S\cap K=\{e\}$, we deduce that $\tilde S\cap\tilde K=\{e\}$.
Since $G=SK$, we have $\tilde G=\tilde S\tilde K$. There is a 
commutative diagram

$$\xymatrix{
\tilde{S}\tilde{K}\ar[d]\ar[r]&SK\ar[d]\\ 
\tilde{K}\ar[r]&K
}$$
where the horizontal arrows are $A$-regular covers and the vertical ones
are projections of trivial left $S$-principal bundles. In other terms,
$\tilde G$ and $G$ are trivial left $S$-principal bundles over $\tilde{K}$ and $K$ resp., and $\tilde G=p^*(G)$ is induced by $G$, where $p:\tilde K\to K$ denotes the restriction of $p:\tilde G\to G$ to $\tilde K$. 

Since $K$ is the connected component of the identity of the real points
of an affine algebraic group, we can identify $K$ with a connected component
of an algebraic subvariety of a Euclidean space ${\mathbb R}^N$.
Since $p:\tilde K\to K$ is a covering map, for each 
$x\in K\subseteq{\mathbb R}^N$, there is a radius $\epsilon_x>0$ such that the
cover is trivial over $B(x,\epsilon_x)\cap K$, and hence we can choose
a continuous section $\sigma_x:B(x,\epsilon_x)\cap K\to\tilde K$ of
$p$. Since $K$ is compact, we can choose a finite trivializing cover
$$K=\bigcup_{i=1}^mK_i,$$
where $K_i=B(x_i,\epsilon_i)\cap K$. Let 
$$A_i=K_i\setminus\bigcup_{j<i}K_j.$$ 
Each set $A_i$ is semi-algebraic, and the sets $A_1,\cdots,A_m$ form 
a partition of $K$. As each $A_i$ is contained in $K_i$,
we obtain by restriction a continuous section 
$\sigma_i:A_i\to\tilde K$ of $p$.
We extend each section $\sigma_i$ by zero outside of $A_i$ and we define a 
section $\sigma:K\to\tilde K$ by $\sigma(x)=\sum_{i=1}^m\sigma_i(x)$. 
We define
$$\tau:G\to\tilde G,\,g=sk\mapsto\tilde s\sigma(k).$$
The map $\tau$ is a Borel section of $p$ and  
is $S$-equivariant: for $s\in S$ and $x\in G$, 
$\tau(sx)=\tilde s\tau(x).$
Let $c:G\times G\to A$ be the corresponding 2-cocycle. That is, 
$c$ is defined by the equation
$$\tau(g)\tau(h)=c(g,h)\tau(gh)\hbox{ for all }g,h\in G.$$
Note that $c$ is a Borel map. It remains to show that $c$ takes only finitely many values in $A$. 

As $G=SK$, taking inverses, 
we also have that $G=KS$. Hence each 
element $g$ of $G$ can uniquely be written as $g=ks$ with $k\in K$ and 
$s\in S$. Let $$\mu:K\times S\to G,\,(k,s)\mapsto ks,$$
denotes the multiplication and let $\pi(\R):G\to G/K$ 
be the canonical projection. 
The quotient space $G/K$ is naturally homeomorphic to $(G_a/K_a)(\R)$,
the connected topological space of $\R$-points of the quotient
variety $G_a/K_a$, and $\pi: G_a\to G_a/K_a$
 is a morphism of 
${\mathbb R}$-varieties \cite[Ch. 2, Thm. 6.8]{BorLin}.
Consider the following morphisms 
$$\theta: S_a\subseteq G_a\to G_a/K_a,$$
between ${\mathbb R}$-varieties \cite[Ch. 2, Prop. 6.12]{BorLin}.
Since $G=SK$ and $S\cap K=1$ the above ${\mathbb R}$-morphism
induces a homeomorphism $S\to
G/K$.  Hence the map 
$$q_S:G\to S,\,sk\mapsto s,$$ 
(the composition of $G\to G/K$ with the
inverse  of $S\to G/K$) is induced by a morphism of ${\mathbb R}$-varieties.
Let ${\bf s}=q_S\circ\mu$. Similarly, the homeomorphism
$$K\to S\backslash G,$$
is obtained by restricting a morphism of ${\mathbb R}$-varieties to topological
connected components.
We define 
$$q_K:G\to K,\,sk\mapsto k,$$ 
and ${\bf k}=q_K\circ\mu$. Let $k\in K$ and 
$s\in S$, there is a unique $s'\in S$ and a unique $k'\in K$ such that
$ks=s'k'$.
By definition of ${\bf s}$ and ${\bf k}$, $s'={\bf s}(k,s)$ and 
$k'={\bf k}(k,s)$, hence
$$ks={\bf s}(k,s){\bf k}(k,s).$$
Let us define
\begin{eqnarray*}f_1:\tilde K\times\tilde K\times\tilde K\times\tilde S\times\tilde S &\to &\tilde G\\
 (u,v,w,p,q)&\mapsto& upvw^{-1}q^{-1}.\end{eqnarray*}
This is a continuous function. Let us define
\begin{eqnarray*}f_2:K\times K\times K\times S\times S &\to &\tilde K\times\tilde K\times\tilde K\times\tilde S\times\tilde S\\
(x,y,z,s,t)&\mapsto &(\sigma(x),\sigma(y),\sigma(z),\tilde s,\tilde t\,).\end{eqnarray*}
This is not a continuous function because $\sigma$ is not necessarily 
continuous, but the restriction of $f_2$ to each semi-algebraic set 
$A_i\times A_j\times A_l\times S\times S,\,1\leq i,j,l\leq m$, is
continuous. 
Let us define
\begin{eqnarray*}f_3:K\times K\times S &\to & K\times K\times K\times S\times S\\
(x,y,s)&\mapsto& (x,y,{\bf k}(x,s)y,s,{\bf s}(x,s)).\end{eqnarray*}
This is an ${\mathbb R}$-morphism, hence it is continuous and the
inverse images 
$$f_3^{-1}(A_i\times A_j\times A_l\times S\times S),\,1\leq i,j,l\leq m,$$
are again semi-algebraic sets. 
Semi-algebraic sets have only a finite number of connected 
components \cite{BR},
say $C_1,\cdots,C_n$ in our case. The sets $C_1,\cdots,C_n$ form a
partition of 
$K\times K\times S$ (since the 
$f_3^{-1}(A_i\times A_j\times A_l\times S\times S)$ do, and the
latter are 
partitioned by their connected components) and the restriction of 
$f=f_1\circ f_2\circ f_3$ to each $C_i$ is continuous.

It is easy to check that $f(x,y,s)=c(x,sy)$ for all $(x,y,s)\in K\times K\times S$
(we defined the functions $f_i$ in order to obtain this equality).
Since $\tau$ is $S$-equivariant and $A$ is central,
$c(sx,ty)=c(x,ty)$ for all $s,t\in S$ and $x,y\in K$.
We conclude that the image of $c$ is equal to the image of $f$.
As $A$ is discrete, it follows that the cardinality of the image of
$f$ 
is bounded by $n$ because $K\times K\times S$ is partitioned by 
the connected sets $C_1,\cdots,C_n$.
\end{proof}
 \begin{coro}\label{linredok} Let $G$ be a virtually connected Lie group with $G^0$ semi-simple (not necessarily linear). Then $G$ has all its Borel cohomology in degree 2 with $\Z$-coefficients bounded.\end{coro}
\begin{proof} We first reduce to the case of the connected component $G^0$
by Proposition \ref{finite index}. Dividing by the center, we reduce then to the case we reduce to the case of $G^0/Z(G^0)$ by using
Corollary \ref{covers}. Now $G^0/Z(G^0):=H$ is a connected linear semi-simple Lie group.
According to \cite{LW} or \cite[Proposition 3.1.6]{Zim}, $H$ is as a Lie group
isomorphic to $G_a(\R)^0$ for some connected real linear algebraic
group $G_a$. It is well known that such a group admits a decomposition
satisfying the assumptions of Theorem \ref{theo: bounded extensions} and we are done.\end{proof}
\begin{rema} In case $G$ is a Lie group of Hermitian type,
it was proved in \cite[Prop. 7.7]{BIW} that the comparison map
$H^2_{Bb}(G,\Z)\to H^2_B(G,\Z)$ is actually an isomorphism.
\end{rema}

\section{Proofs of the main results}\label{sect: Proofs of the main results}

\subsection{Proof of Theorem \ref{theo: higher classes}} We first need a lemma.
\begin{lemma}\label{res} Let $G$ be a connected Lie group an $\sqrt G$ its radical.
Then the restriction map $H^*(BG,\R)\to H^*(B\sqrt G,\R)$ is surjective.
\end{lemma}
\begin{proof}
According to Proposition \ref{prop: onto}, the map 
$$H^2(BG,\R)\to H^2(B\sqrt G,\R),$$ 
is surjective. Because
$\sqrt G$ is homotopy equivalent to a torus, the cohomology $H^*(B\sqrt G,\R)$ is
a polynomial algebra on 2-dimensional classes. Therefore,
since the restriction map $H^*(BG,\R)\to H^*(B\sqrt G,\R)$ is a map of algebras, it is surjective.
\end{proof}
\begin{proof} We prove Theorem \ref{theo: higher classes}. Let $I^*$ denote the image of the map $H^*(B(G/\sqrt G),\R)\to H^*(BG,\R)$.
According to Lemma \ref{res}, the restriction map
$H^*(BG,\R)\to H^*(B\sqrt G,\R)$ is surjective and therefore the fiber in the fibration
$B\sqrt G\to BG\to BG/\sqrt G$ is totally non-homologous
to zero relative to $\R$.
We conclude from Borel's Theorem \cite[Thm. 14.2]{BorTop})
that $H^*(BG,\R)/\left< I^{+}\right>$ maps isomorphically onto $H^*(B\sqrt G,\R)$,
where $\left< I^{+}\right>$ stands for the ideal generated by the classes 
$I^+<I^*$
of positive degree. Using the well-known fact that for a connected Lie group $L$, $H^*(BL,\R)$
is a polynomial algebra on even dimensional generators
(see \cite[Thm. 18.1]{BorTop})
we now prove our claim
on the generation of $H^*(BG,\R)$ by checking it for generators
in even degrees. Thus, assume
that $x\in H^{2d}(BG,\R)$ is a generator with $d>1$. Then $x$ maps 
to an element $y\in
H^{2d}(B\sqrt G,\R)$, which is a product of $d$ 2-dimensional classes $y_i$.
By choosing counter images $x_i$ in $H^2(BG,\R)$ for the $y_i$'s, 
their product will be
an element $z\in H^{2d}(BG,\R)$ with $x-z\in\left< I^{+}\right>$, say
$x-z=\sum a_j b_j$ with $a_j\in I^{+}$ and $b_j\in H^{<2d}(BG,\R)$.
Thus $x=z+\sum a_j b_j$ lies in the subring generated
by $H^2(BG,\R)$ together with $\left< I^+\right>$, proving the theorem. \end{proof}

\subsection{Proof of Theorem \ref{theo: main}}
\subsubsection{If the class defined by the universal cover is bounded then the fundamental group is undistorted}
\begin{proof} We prove that (\ref{Uni}) implies (\ref{Dis}). Let $A\cong \pi_1(G)$ denotes the fundamental group of $G$. We have a topological $A$-extension:
	$$0\to A\to\tilde{G}\to G\to 1,$$
where $\tilde G$ denotes the universal cover of $G$. Let $\tau$ be a Borel section of the projection $\tilde{G}\to G$. By hypothesis the class
$[c_{\tau}]\in H^2_B(G,A)$ is bounded. According to Lemma \ref{lemm: cocycle defined by a section}, there exists a Borel section $\sigma$ of $\tilde{G}\to G$
such that $c_{\sigma}$ is bounded. We may modify $\sigma$ at the identity so that $\sigma(e)=e$. It is easy to check that the new cocycle $c_\sigma$ associated to the
modified section $\sigma$ is still bounded. Proposition \ref{prop: undistorted} implies that $A$ is undistorted in $\tilde G$.	
\end{proof}
\subsubsection{If the fundamental group is undistorted then the radical is linear}
\begin{proof} We prove that (\ref{Dis}) implies (\ref{RadLin}). Let $\sqrt G$ be the radical of $G$. According to Proposition \ref{prop: general equivalences}, the non-linearity of $\sqrt G$ implies the existence of a maximal compact subgroup $K$ of $\sqrt G$ such that
$[\sqrt G,\sqrt G]\cap K\not=\{e\}$.
Hence, Proposition \ref{prop: distorted} applies, and we deduce that the fundamental group is distorted.
\end{proof}
\subsubsection{If the $2$-dimensional cohomology is bounded then the radical is linear}
\begin{proof} We show that (\ref{H2b}) implies (\ref{RadLin}). Let us assume that $\sqrt G$ is not linear. We will show that $H^2_B(G,\Z)$ contains a class with no Borel bounded representative.
According to Proposition \ref{prop: general equivalences}, the non-\-linea\-rity of $\sqrt G$ implies the existence of a maximal compact subgroup $K$ of $\sqrt G$ such that
$[\sqrt G,\sqrt G]\cap K\not=\{e\}$.
Hence, Proposition \ref{prop: distorted} applies, and we deduce the existence of
a distorted topological $\Z$-extension 
$$0\to\Z\to E\to G\to 1.$$ 
Consider a Borel cocycle $c$, such that $[c]\in H^2_B(G,\Z)$ is defined by the extension $E$ of $G$.
We claim that $c$ is not bounded. Assume this is not the case.
Lemma \ref{lemm: cocycle defined by a section} shows that $c=c_{\sigma}$, where $\sigma$ is some Borel section of $E\to G$ (and we may modify $\sigma$ at the identity so that $\sigma(e)=e$, keeping the modified cocycle $c_\sigma$  bounded). Proposition \ref{prop: undistorted} implies that the extension $0\to\Z\to E\to G\to 1$ is undistorted, contradicting the assumption. This proves that $c$ is unbounded.
\end{proof}
\subsubsection{If the radical is linear then the $2$-dimensional cohomology is bounded}
We first prove the following lemma. 
\begin{lemma}\label{fix} Let $G$ be a connected Lie group. Assume that the radical $\sqrt G$ of $G$ is linear. Then there is
a closed contractible subgroup $V<\sqrt G$ which is normal in $G$ such that 
$\sqrt G/V$ is a torus $T$ and the covering space $\xi$ induced from the universal cover $\tilde{L}\to G/\sqrt G$ via the canonical projection $G/V\to G/\sqrt G$ is of the form 
$\xi: T\times\tilde{L}\to G/V$.\end{lemma}
\begin{proof}
Let $H$ be the covering space of $G$ induced from the universal cover $\tilde{L}$ of $G/\sqrt G$. Then, according to (\ref{pullback exact sequence}) of Lemma \ref{lemm: pullback} and Lemma \ref{gensplitting}, $H= \sqrt G\rtimes \tilde{L}$, with $\tilde{L}$ semi-simple. Let $W$ be the linearizer of $\tilde{L}$; it is discrete and central in
$\tilde{L}$. Also $W$ is central in $H$ because $W$ lies in the kernel of the adjoint representation of $H$. So we can form the quotient 
$$Q:=\sqrt G\rtimes (\tilde{L}/W),$$ 
which is a linear Lie group by Theorem \ref{OV}. According to Hochschild \cite[XVIII.4, Thm. 4.3]{HocTheStr}, the group $Q$ contains 
a contractible normal closed solvable subgroup $X$ with $Q/X$ linear reductive. Since $X<\sqrt G$ we can think of $X$
as a subgroup of $H=\sqrt G\rtimes\tilde{L}$, and this subgroup is normal, since $W$ is central. Under the covering
projection $H\to G$, this subgroup $X$ maps isomorphically onto a closed normal subgroup $V<G$, which lies
in the radical of $G$. Since $Q/X$ is linear reductive, $\sqrt G/X$ is a central torus $T$ (see \cite[XVIII.4, Thm. 4.4]{HocTheStr}). Let $\xi$ be the covering induced from the universal cover $\tilde{L}\to G/\sqrt G$ via the canonical projection $G/V\to G/\sqrt G$. According to (\ref{pullback exact sequence}) of Lemma \ref{lemm: pullback}, the total space $E$ of $\xi$ sits in the short exact sequence
$$0\to T\to E\to \tilde{L}\to 1.$$
As $\tilde{L}$ is simply-connected, and as $T$ is the radical of $E$, the sequence splits (again by Lemma \ref{gensplitting}). Since $T$ is central, we deduce that $E$ is a direct product $E=T\times\tilde{L}$.\end{proof}

\begin{proof} We prove that (\ref{RadLin}) implies (\ref{H2b}).
By Lemma \ref{fix}, $G$ contains a normal contractible
subgroup $V$, hence $G\to G/V$ is a homotopy equivalence, and thus $H^*_B(G/V,\Z)\to
H^*_B(G,\Z)$ an isomorphism by Lemma \ref{hinv}. It suffices therefore to prove (\ref{H2b}) for $G/V$. In the
notation of Lemma \ref{fix}, we have covering spaces $\R^n\times \tilde{L}\to T\times
\tilde{L}\to G/V$. Note that $\R^n\times\tilde{L}$ is the universal cover of $G/V$. By
Proposition \ref{prop: universal cover}, it suffices to show that the cocycle defining this universal cover can
be chosen
bounded. For this, we appeal to Lemma \ref{lemm: composition}, with $X=\R^n\times\tilde{L}$, $Y= T\times\tilde{L}$ and $Z=G/V$. The map $p:X\to Y$ is the product of $\R^n\to T$ with
${\rm Id}:\tilde{L}\to\tilde{L}$ and we conclude from Lemma \ref{lemm: Kunneth}
and Corollary \ref{compact}
that the Borel cocycle
associated with $p$ can be chosen to be bounded. For $q:Y\to Z$ we observe that $q$ is obtained 
by pulling back $\tilde{L}\to G/\sqrt G$ over $G/V$. According to (\ref{pullback section}) of Lemma \ref{lemm: pullback},
it suffices
to prove that the universal cover $\tilde{L}\to G/\sqrt G$ has a Borel section
$\sigma$ such that the associated cocycle $c_{\sigma}$ is bounded. 
But this case has already been dealt with in Corollary \ref{linredok}
and we are done.
\end{proof}

\subsubsection{If the $2$-dimensional cohomology is bounded then all the cohomology is bounded}
\begin{proof} We prove that (\ref{H2b}) implies (\ref{Hb}).
Since for every $d$, $H^d(BG,\Z)$ is a finitely generated abelian group,
we see that the $\Z$-subalgebra $A^*\subseteq H^*(BG,\R)$ generated
by $H^2(BG,\Z)$ together with the image $J^*$ of the inflation map
$H^*(B(G/\sqrt G),\Z)\to H^*(BG,\Z)$ satisfies, in view of
Theorem \ref{theo: higher classes}, $A^*\otimes \R=
H^*(BG,\R)$. Thus $A^d\subseteq H^d(BG,\Z)$ has finite index for
every $d$. To show that $H^d_{Bb}(G,\Z)\to H^d_B(G,\Z)$ is 
onto it suffices therefore by Lemma \ref{modn} to show that the elements
of $A^d$ all lie in the image of $H^d_{Bb}(G,\Z)\to H^d_B(G,\Z)$.
From the definition of $A^*$ it suffices therefore to show that:
\begin{itemize}
\item[(a)] The map $H^2_{Bb}(G,\Z)\to H^2_B(G,\Z)$ is onto.
\item[(b)] The map $H^*_{Bb}(G/\sqrt G,\Z)\to H^*_B(G/\sqrt G,\Z)$ is onto.
\end{itemize}
Now (a) holds, because of our assumption (2). To prove (b), we consider
the Lie group $L$ obtained by dividing $G/\sqrt G$ by its
center. By Corollary \ref{covers} it suffices to show
that $H^*_{Bb}(L,\Z)\to H^*_B(L,\Z)$ is onto. But this
follows from the second statement in Proposition \ref{prop: michelle},
by observing that $L$ is as a Lie group isomorphic to
the topological connected component of the identity of a Lie group $G_a(\R)$
for some real linear algebraic group $G_a$ (cf. \cite{LW} or \cite{Zim}).
\end{proof}
\subsubsection{Remaining implications}
That (\ref{Hb}) implies (\ref{H2b}) is a tautology.
That (\ref{H2b}) implies (\ref{Uni}) follows from Proposition \ref{prop: universal cover}.
That (\ref{RadLin}) is equivalent to (\ref{Com}) follows from Proposition
\ref{prop: general equivalences}.

\subsection{Proof of Theorem \ref{theo: vir}}
\begin{proof}
In view of Proposition \ref{finite index} we may assume that
$G$ is connected. Hence Theorem \ref{theo: vir} follows from
Theorem \ref{theo: main}.
\end{proof}
\subsection{Proof of Theorem \ref{theo: generalization of Goldman}}
\subsubsection{Goldman's construction}
Recall that if $G$ is a connected Lie group and if $P$ is a $G$-principal  bundle over a $CW$-complex $B$, there is a characteristic class:
$$o_2(P)\in H^2(B,\pi_1(G)),$$
defined by obstruction theory. (More precisely, let $e^{(2)}$ be a $2$-cell of the $2$-skeleton $B^{(2)}$ of $B$. There is, up to homotopy, a unique way
of getting $P|_{B^{(2)}}$ by gluing $P|_{B^{(2)}\setminus interior(e^{(2)})}$
with $e^{(2)}\times G$. The gluing is determined by a map 
$\partial e^{(2)}\to G$. The homotopy classes $c(e^{(2)})\in\pi_1(G)$ of these maps 
define a cellular cocycle whose class is $o_2(P)$.) 

\begin{prop}\emph{(Compare with Goldman \cite{Gol}.)}\label{prop: Goldman's construction}
Let $G$ be a connected solvable Lie group. If $\pi_1([G,G])\neq 0$, then there exists a class of degree $2$ in the image of the natural map $H^*(BG,\R)\to H^*(BG^{\delta},\R)$ which can't be represented by a bounded cocycle.  
\end{prop}

\begin{proof} Let $p:\tilde G\to G$ be the universal cover of $G$. The commutator subgroup $[\tilde G,\tilde G]$ is closed connected and normal in $\tilde G$. Lemma \ref{lemm: pi1} implies that $[\tilde G,\tilde G]$ is simply-connected. Obviously $p([\tilde G,\tilde G])=[G,G]$, hence the restriction of $p$ to $[\tilde G,\tilde G]$ is the universal cover of the subgroup $[G,G]$ of $G$. (Here the topology on the group $[G,G]$ is the quotient topology induced by the restriction of $p$ to $[\tilde G,\tilde G]$.)
Hence we obtain two exact sequences with commutative squares:
$$\xymatrix{
0\ar[r]&\pi_1(G)\ar[r]&\tilde G\ar[r]^p&G\ar[r]&1\\
0\ar[r]&\pi_1([G,G])\ar[r]\ar[u]&[\tilde G,\tilde G]\ar[r]\ar[u]&[G,G]\ar[r]\ar[u]&1.
}$$
Each vertical arrow is an inclusion (but $[G,G]\subset G$ is not necessarily closed). A connected solvable Lie group has the homotopy type of its maximal
compact subgroup which is a torus. Hence, by hypothesis, there exists 
an element $z\in \pi_1([G,G])\subset Ker(p)$ of infinite order. Let $a_1,\dots, a_g, b_1,\dots, b_g$ be elements of $\tilde G$ such that
$$z=\prod_{k=1}^g[a_k,b_k].$$
Let $\Sigma_g$ denote the orientable closed surface of genus $g$ and let
$$\Gamma=\pi_1(\Sigma_g)=(A_1,\dots, A_g, B_1,\dots, B_g: \prod_{k=1}^g[A_k,B_k]),$$ be its fundamental group with the usual presentation. Let $h:\Gamma\to G$
be the homomorphism defined by the conditions $h(A_k)=p(a_k),\, h(B_k)=p(b_k)$.
Let $P_h=\Gamma\backslash(\tilde\Sigma_g\times G)$ be the flat principal $G$-bundle
over $\Sigma_g$ with $h$ as holonomy. In other words, the bundle $P_h$ is the pullback of the universal principal $G$-bundle $EG$ by
the map $Bh:B\Gamma\to BG$.
By construction:
$$o_2(P_h)([\Sigma_g])=z.$$
As $z$ has infinite order in the finitely generated abelian group $\pi_1(G)$, we can choose a homomorphism $q:\pi_1(G)\to\R$  between coefficients such that
$q(z)=1$. Let $q_*:H^2(BG,\pi_1(G))\to H^2(BG,\R)$ be the induced map.
We have $q_*(o_2(P_h))([\Sigma_g])=1$.
The map $Bh$ factors through $BG^{\delta}$:
$$\Sigma_g\simeq B\Gamma\rightarrow BG^{\delta}\rightarrow BG.$$
It induces the commutative diagram: 
$$\xymatrix{
P_h\ar[d]\ar[r]&EG^{\delta}\ar[d]\ar[r]&EG\ar[d]\\
\Sigma_{g}\ar[r]&BG^{\delta}\ar[r]&BG,}$$
where vertical arrows are projections of principal $G$-bundles.
It also induces the commutative diagram:
$$\xymatrix{
H^2(BG,\pi_1(G))\ar[d]\ar[r]&H^2(BG^{\delta},\pi_1(G))\ar[d]
\ar[r]&H^2(\Sigma_g,\pi_1(G))\ar[d]\\
H^2(BG,\R)\ar[r]&H^2(BG^{\delta},\R)
\ar[r]&H^2(\Sigma_g,\R),}$$
in which characteristic classes are preserved and where vertical arrows are induced by the coefficients homomorphism $q$.
The image $x^{\delta}\in
H^2(BG^{\delta},\R)$ of $x:= q_*(o_2(EG))\in H^2(BG,\R)$ is not equal to zero 
because it is sent to $q_*o_2(P_h)$ and we have proved that this class is non-zero by
evaluating it on $[\Sigma_g]$. The class $x^{\delta}\in H^2(BG^{\delta},\R)\cong H^2(G^{\delta},\R)$ can't be bounded because a solvable Lie group is amenable
also with respect to its discrete topology \cite[Volume 1, Chapter IV, Paragraph 17,
(17.14)]{HR}. Hence its  real bounded cohomology vanishes in positive dimensions.
\end{proof}

\begin{proof} We prove Theorem \ref{theo: generalization of Goldman}.
We apply Proposition \ref{prop: Goldman's construction} to the group $\sqrt G$. We conclude that there exists 
a class $x\in H^2(B\sqrt G ,\R)$ with $x^{\delta}\in H^2(B(\sqrt G)^{\delta},\R)$
not bounded. According to Proposition \ref{prop: onto}, the map
$H^2(BG,\R)\to H^2(B\sqrt G,\R)$ is onto. Let $y\in H^2(BG,\R)$ be a counter image of
$x$. The commutativity of the diagram:
$$\xymatrix{
H^2(BG,\R)\ar[d]\ar[r]&H^2(B\sqrt G,\R)\ar[d]\\
H^2(BG^{\delta},\R)\ar[r]&H^2(B(\sqrt G)^{\delta},\R),}$$
implies that $y^{\delta}$ is unbounded.
\end{proof}

\subsection{Proof of Theorem \ref{theo: Goldman}}
\subsubsection{Condition \ref{commutator} implies Condition \ref{virtually trivial}}
\begin{proof} That (\ref{commutator}) implies (\ref{virtually trivial}) is one of  of the implications of Goldman's result in \cite{Gol}, so we don't prove it but
we show with an example that the finiteness assumption on the complexes in (\ref{virtually trivial}), which is not explicitly stated as an hypothesis in \cite{Gol} (but used in the proof), is needed. 
Let $G=SO(2,\R)$. 
The bundle $EG^{\delta}\to BG^{\delta}$ induced by the map $BG^{\delta}\to BG$ is a flat principal $G$-bundle.
Although $[G,G]=1$, this bundle
is not virtually trivial because  the image $e^{\delta}\in H^2(BG^{\delta},\Z)$ of the Euler class $e\in H^2(BG,\Z)$
has infinite order \cite[Cor. 1]{MilOnTheHom}.
\end{proof} 
\subsubsection{Condition \ref{virtually trivial} implies Condition \ref{zero}}
\begin{proof} 
Let $c(EG)\in H^*(BG,\Z)$ be a characteristic class. Let $P$ be a flat principal
$G$-bundle over a finite complex $X$. By hypothesis there exists a finite
cover $X_d\to X$, say of degree $d$, such that the pullback $P_d$ of $P$ to $X_d$ is trivial. Hence, according to Lemma \ref{lemm: transfer for covers}, we have $d\cdot c(P)=0$. We apply (\ref{restriction is torsion}) of Proposition \ref{prop: flat bundles} and conclude that $c(EG)$ is sent to zero in  $H^*(BG^{\delta},\R)$. As  $H^*(BG,\Z)\otimes\R=H^*(BG,\R)$  this shows that
the map $H^*(BG,\R)\to H^*(BG^{\delta},\R)$ is zero.	
\end{proof}
\subsubsection{Remaining implications}
That (\ref{zero}) implies (\ref{bounded}) is a tautology. That (\ref{bounded}) implies (\ref{commutator}) follows from Theorem \ref{theo: generalization of Goldman}.

\end{document}